\documentclass{article}

\usepackage{arxiv}

\usepackage[utf8]{inputenc} 
\usepackage[T1]{fontenc}    
\usepackage{hyperref}       
\usepackage{url}            
\usepackage{booktabs}       
\usepackage{amsfonts}       
\usepackage{nicefrac}       
\usepackage{microtype}      
\usepackage{lipsum}

\usepackage{amssymb}
\usepackage{bm}
\usepackage{amsmath}
\usepackage{lineno}
\usepackage{graphicx}
\usepackage{subcaption}

\newcommand{\B}{\bm{B}}
\renewcommand{\v}{\bm{v}}
\newcommand{\m}{\bm{m}}
\renewcommand{\k}{\bm{k}}
\newcommand{\A}{\bm{A}}
\newcommand{\I}{\bm{I}}
\renewcommand{\u}{\bm{u}}
\newcommand{\PHI}{\bm{\Phi}}
\newcommand{\LAMBDA}{\bm{\Lambda}}
\renewcommand{\c}{\bm{c}}
\newcommand{\x}{\bm{x}}

\title{A General Shock Waveform and Characterisation Method}

\author{
  Yinzhong Yan\\
School of Mechanical, Aerospace and Civil Engineering\\
 The University of Manchester\\
  Manchester, M13 9PL, United Kingdom\\
  \texttt{Yinzhong.Yan@manchester.ac.uk} \\
   \And
 Q.M. Li\thanks{Corresponding author.}  \\
School of Mechanical, Aerospace and Civil Engineering\\
The University of Manchester\\
Manchester, M13 9PL, United Kingdom\\
\texttt{Qingming.Li@manchester.ac.uk} \\
}

\begin{document}
\maketitle

\begin{abstract}
A shock waveform is proposed based on the mechanical mechanism of shock generation in a structure.
The parameters in the shock waveform have clear mechanical meanings about the generation and development of the shock.
A shock signal processing method is proposed and applied to represent a pyroshock or ballistic shock signal in both temporal and frequency domain using finite terms of the shock waveform components.
It is found that complexity and dominant shock distance of a shock can be described quantitatively by the number of waveform components ($\eta_{90\%}$) and a normalised parameter ($\kappa$), respectively.
Pyroshock and mechanical shocks in different categories are analysed using the proposed shock waveform to demonstrate its importance and generality in shock representation.
\end{abstract}


\section{Introduction}
Shock environment has to be considered for the design of equipment and components used in aerospace, transportation, defence and package engineering.
Shock signals are normally recorded in the form of acceleration-time-history data, which is not straightforward to be linked to its severity and its effect on equipment components\cite[p.385]{ECSS2015}.
It is often necessary to get more specific mechanical information from shock signal using signal processing tools to understand the features of the shock and the associated phenomena.

Shock signal can be decomposed into a combination of basis signals.
The most well-known decomposition method is Fourier transform, with which the frequency content of a measured signal can be identified by representing the signal with a series of harmonic waveforms.
NASA proposed a decomposition method and an associated algorithm for shock signal using wavelet\cite{ferebee2008}, with intention to replace the shock response spectrum (SRS) for shock specification.
Prony decomposition method is usually adopted to express the shock as the superposition of exponentially damped harmonic waves in recent research study\cite{monti2017,hwang2016,yan2019}.
The Prony decomposition has advantage to describe near-field shocks characterised by the simultaneous peak location and sharp initial rise of its components.
ESA provided an advanced scheme for Prony decomposition\cite{ECSS2015,bernaudin2005}, where the advanced Prony mode is given by the standard Prony mode convoluted with a basic Gaussian pulse.
This scheme has advantage for the far-field shock environment characterised by the different peak time and the gradual initial rise of each advanced Prony mode. 
However, ESA did not give explicit mathematical expression for this scheme and its decomposition algorithm is not available.
Similarly, the decomposed components of shock signals by empirical mode decomposition\cite{wang2011transient}, discrete wavelet decomposition or wavelet package decomposition methods\cite{singh2018rolling} are discrete signals without mathematical expressions. 

Other basis waveforms can be used to synthesise mechanical shock environment for design and testing purposes.
The WAVSYN waveform was initially proposed for the simulation of earthquakes, and was also found effective in shock synthesise\cite{yang1970safeguard,yang1972development}.
A waveform function similar to the response of a single-degree-of-freedom system to a Dirac impulse excitation was proposed in Ref.\cite{kern1984}.
The ZERD waveform function was proposed to simulate shock signals with zero residual displacements, which has been well adopted with a shaker\cite{fisher1977,smallwood1985}.
A combination of simple impulses, sine waves, damped sine waves and modified Morlet wavelet has been used to synthesise shock time histories\cite{brake2011}.
These methods are mainly for the synthesis of shock signal corresponding to a given SRS, rather than for the decomposition of a measured shock time history signal.

This study analyses the features of the generation mechanism of pyroshock and ballistic shock, and proposes a characteristic shock waveform to represent shock signals.
This new shock waveform is determined by a group of parameters, i.e., amplitude, frequency, damping ratio, initial time, peak time and phase, which have explicit mechanics meaning.
By decomposing the shock signal into the sum of limited terms of shock waveform components, mathematical expression of a shock signal can be explicitly given, which can support the understanding of shock effects.

\section{Shock Waveform Expression from Shock Response}

\subsection{Shock Waveform}

The shock waveform is proposed in the form of Eq.(\ref{inital_waveform}) based on the mechanical mechanism described in \ref{mechanical_background}.
\begin{equation}\label{inital_waveform}
W(t)=t^n e^{-\omega \zeta t+i (\omega t+ \varphi)} H(t)
\end{equation}
where $n$ is the order, $\omega$ is the angular frequency, $\zeta$ is the damping ratio, $i$ is the imaginary unit, $\varphi$ is the phase, and the $H(t)$ is the Heaviside step function defined by
\begin{equation}\label{heaviside}
H(t)=
\begin{cases}
0,\quad t<0\\
1,\quad t\geq0
\end{cases}
\end{equation}
Step function is included because the magnitude of shock signal is zero before the start of shock event.

If $W(t)$ is considered as a modulated harmonic signal, the envelope of the shock waveform $W(t)$ is
\begin{equation}
M(t)=t^n e^{-\omega \zeta t} H(t)
\end{equation}
By calculating the stationary point of $M(t)$ in $t>0$, the time $\tau$ associated with the maximum of $M(t)$ is determined by
\begin{flalign}
&& M'(\tau) &=n \tau^{n-1} e^{-\zeta \omega \tau }-\zeta  \omega  \tau^n e^{-\zeta \omega \tau}=0 & \\
\text{or} && \tau &=\frac{n}{\omega \zeta} & \label{introduce_tau}
\end{flalign}
The shock waveform proposed in Eq.(\ref{inital_waveform}) can be normalized in terms of the maximum value of its envelope, i.e.
\begin{align}
\bar{W}(t) & =\frac{W(t)}{M(\tau)}\\
&=\frac{M(t)}{M(\tau)}e^{i (\omega t + \varphi)} \label{normalization}\\
& =t^{\zeta \omega \tau} \tau^{-\zeta \omega \tau} e^{\zeta \omega (\tau-t)+i (\omega t+\varphi)} H(t)
\end{align}
which can be generalised as
\begin{equation}\label{shock_waveform}
w(t)=A t^{\zeta \omega \tau} \tau^{-\zeta \omega \tau} e^{\zeta \omega (\tau-t)+i (\omega t+\varphi)} H(t) 
\end{equation}
where $A$ is the amplitude of the waveform.

With the relationship shown in Eq.(\ref{introduce_tau}), the parameter $n$, which has no direct mechanical meaning in Eq.(\ref{inital_waveform}), is now replaced by the parameter $\tau$ representing the peak time of a single shock waveform.
$\tau$ is an important parameter when analysing the temporal structure of the shock waveform.

Considering a practical shock signal, the initial time for each shock waveform component may be different.
Therefore, in the superposition process, the initial time for each shock waveform should be taken into consideration.
Then, a shock signal $r(t)$ can be generally expressed as the superposition of a series of shock waveform components, i.e.
\begin{equation}\label{decomposition}
r(t)=\sum_{i=1} T_i (w_i(t))
\end{equation}
where $T_i$ is a temporal translation transform for the consideration of the initial time for each waveform component, i.e.,
\begin{equation}\label{translation_transformation}
T_i (w_i(t,\tau, A, \omega, \zeta, \varphi))=w_i(t-\mathring{t}_i,\tau-\mathring{t}_i,A, \omega, \zeta, \varphi)
\end{equation}
and $\mathring{t}_i$ is the initial time of the $i$th shock waveform.

\subsection{Mechanical Significance of Waveform Parameters}

To better illustrate the shock waveform, the mechanical significances of waveform parameters in Eq.(\ref{shock_waveform}) are discussed below.

\paragraph{Frequency $\omega$}
\begin{figure}
	\centering
	\includegraphics[width=0.5\linewidth]{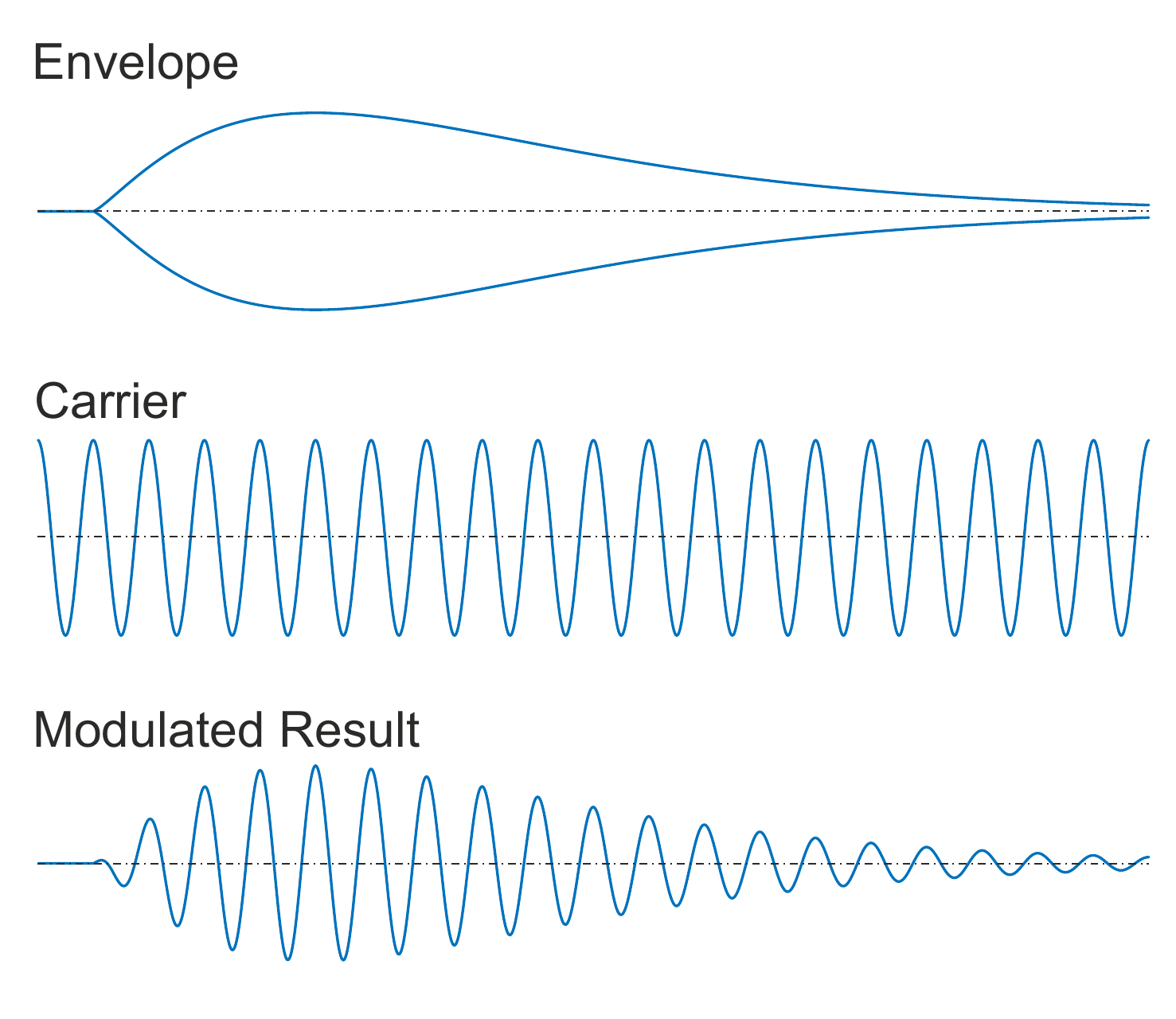}
	\caption{Amplitude modulation of shock waveform}
	\label{diagram_modulation}
\end{figure}
The shock waveform can be considered as an amplitude modulated harmonic wave.
Here $\omega$ is actually the frequency of carrier wave, as depicted in Fig.\ref{diagram_modulation}.
Such a shock waveform will cause resonance for those modes whose natural frequencies are closed to $\omega$.
Therefore, $\omega$ is one of the key parameters of the shock waveform.

\paragraph{Amplitude $A$}

\begin{figure}
	\centering
	\includegraphics[width=0.5\linewidth]{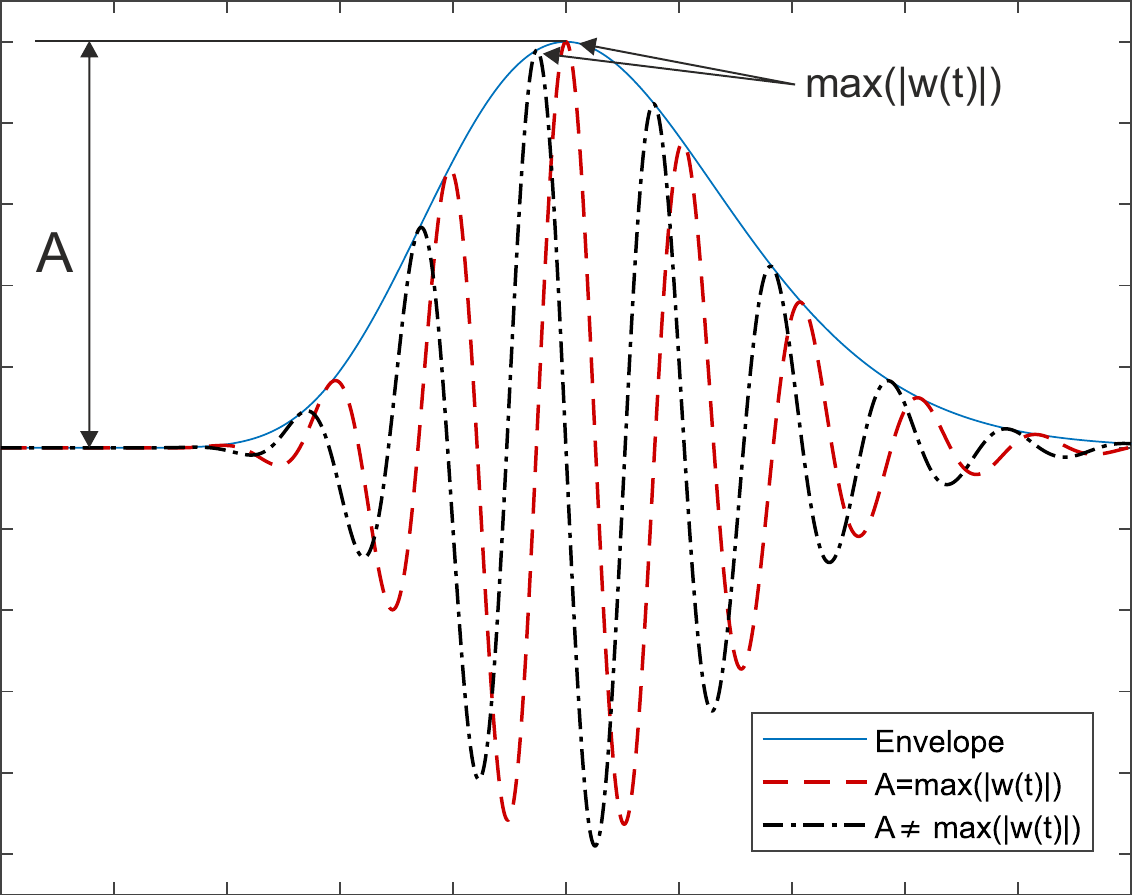}
	\caption{The difference between amplitude $A$ and $max(|w(t)|)$}
	\label{effect_a_varphi}
\end{figure}

This is another key parameter of the shock waveform.
The amplitude $A$ is the maximum of envelope curve, rather than the absolute maximum of a shock waveform $max|w(t)|$.
As shown in Fig.\ref{effect_a_varphi}, usually the difference between maximum value of envelope and shock waveform are small.
The condition of $A=max|w(t)|$ is that absolute maximum of carrier wave meets the maximum of its envelope, which is equivalent to
\begin{equation}\label{maximum_condition}
\varphi=-(\tau \omega-\pi \lfloor \frac{\tau \omega}{\pi} \rfloor)
\end{equation}
where $\lfloor x \rfloor$ is the floor function with input $x$.

\paragraph{Damping ratio $\zeta$}
\begin{figure}
	\centering
	\includegraphics[width=0.5\linewidth]{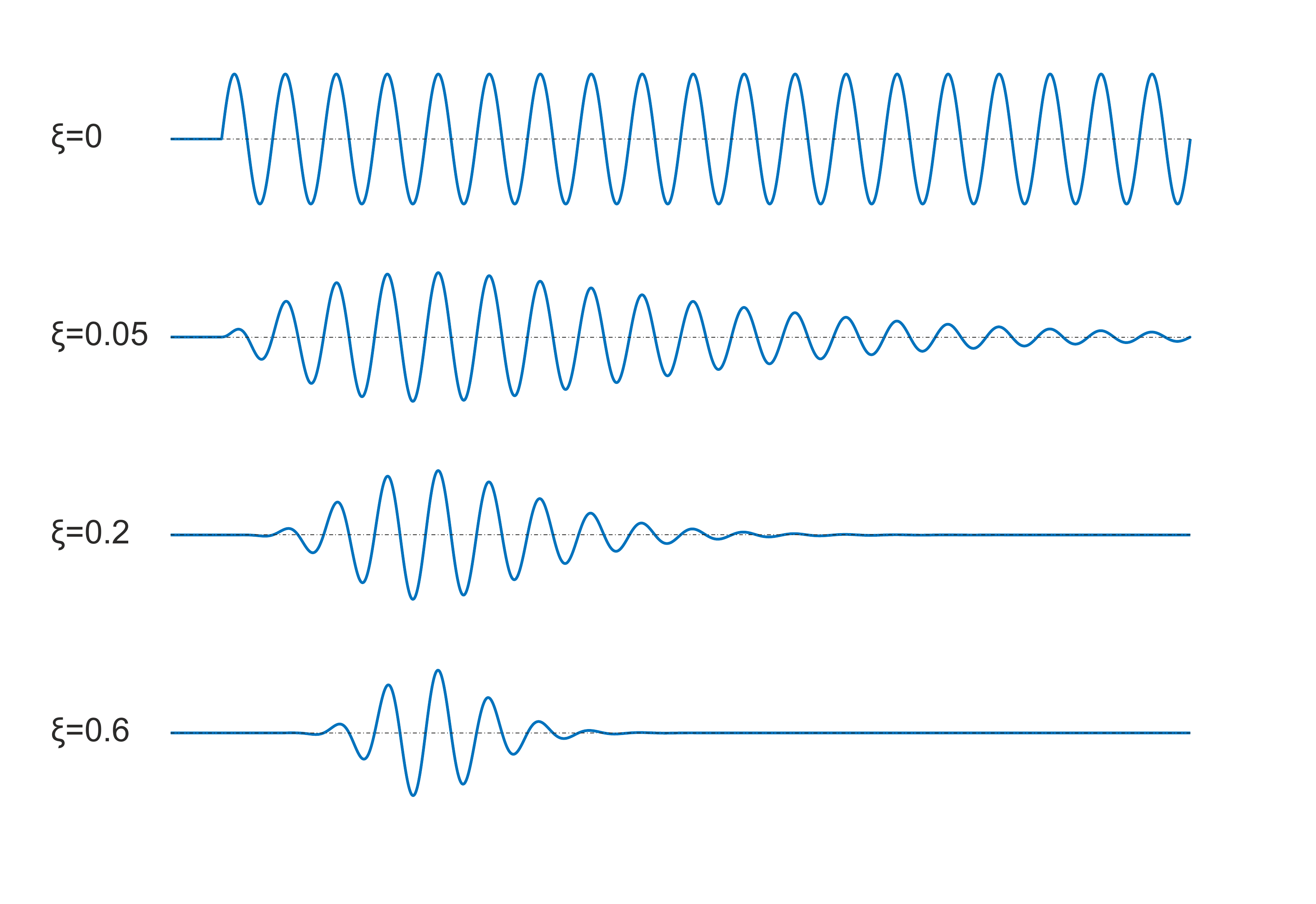}
	\caption{Influences of damping ratio $\zeta$ on duration of shock waveform}
	\label{diagram_damping}
\end{figure}
This parameter mainly influences the duration of shock waveform.
Smaller damping ratio leads to longer shock waveform duration, while larger damping ratio results in shorter duration, as illustrated in Fig.\ref{diagram_damping}.
In the extreme case when $\zeta=0$, shock waveform is a harmonic wave.

\paragraph{Peak Time $\tau$}
\begin{figure}
	\centering
	\includegraphics[width=0.6\linewidth]{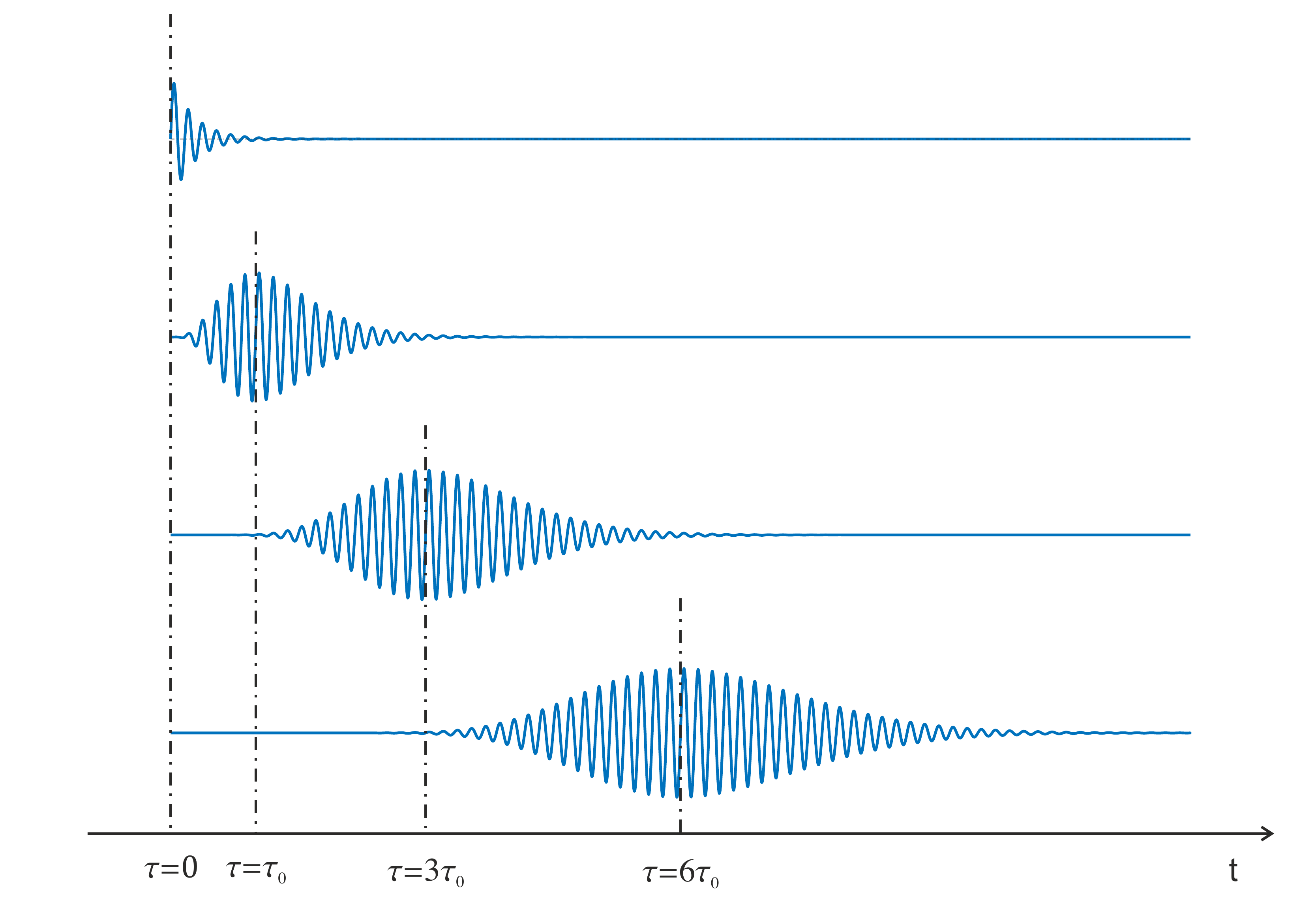}
	\caption{Influences of peak time $\tau$, here $\tau_0=1$ ms}
	\label{diagram_tau}
\end{figure}
As shown in Fig.\ref{diagram_tau}, $\tau$ defines the temporal structure of each waveform component, i.e., the duration between the initial time and the peak time.
If peak time aligns with initial time, i.e., $\tau=0$, shock waveform represent a Prony mode (ordinary damped harmonic wave).
With the increase of $\tau$, the temporal structure of the waveform changes from Prony mode to the advanced Prony mode.
When $\tau$ is large enough, the shape of shock waveform approaches symmetric wavelet.
Their relationship will be further discussed in Section \ref{section_waveform_comparision}.

The physical meaning of peak time $\tau$ can be used to evaluate dominant shock distance, i.e., near-field, mid-field and far-field.
The smaller the peak time, the closer the shock measurement point to the shock source.
Considering the effect of frequency of carrier wave in time domain, it is proposed here that the dominant shock distance can be indicated by a parameter normalised by the period of the carrier harmonic wave, i.e.
\begin{equation}
\kappa=\frac{\tau \omega}{2 \pi}=\frac{\tau}{T}=\frac{\tau c_0}{T c_0}=\dfrac{d}{\lambda}
\end{equation}
where $T$ is the period of the carrier harmonic wave; $d$ is the distance between the excitation source and the measurement point; $c_0$ and $\lambda$ are the phase velocity and the wavelength of the wave, respectively.
Therefore, the parameter $\kappa$ is equivalent to the distance between the excitation source and the measurement point normalised by the wavelength of the wave.

\begin{table}
	\centering
	\caption{Shock environment classification based on $\kappa$}
	\label{environment_category}
	\begin{tabular}{|c|c|}
		\hline 
		Category & $\kappa$ Range \\ 
		\hline 
		Near-field shock & $0 \leqslant \kappa < 1$ \\ 
		\hline 
		Mid-field shock & $1 \leqslant \kappa < 10$ \\ 
		\hline 
		Far-field shock & $\kappa \geqslant 10$ \\ 
		\hline 
	\end{tabular}
\end{table} 

In this study, a quantitative classification of the dominant shock distance is given in Table \ref{environment_category}, based on the order of the magnitude of $\kappa$.
With this classification, the waveforms of the near-field shock are similar to the damped sine waves, while the waveforms of the far-field shocks are close to the advanced Prony mode (or wavelet if $\zeta$ is high).
Similar to the definition in Ref.\cite{ECSS2015}, the transitional shocks with characters between near- and far-field shocks are classified as mid-field shock.
Further analysis and discussion about $\kappa$ are given in Section \ref{section_analysis}.

\paragraph{Initial Time $\mathring{t}$}

This parameter is necessary for the reconstruction of shock signal by the superposition of shock waveform component.
The proposed shock waveform $w(t)$ has $t=0$ as the initial time, while it is possible that a shock signal is composed of shock waveform components having different non-zero initial times.
Therefore, the simple translation transform in time domain can move the waveform component from $0$ to $\mathring{t}$.
The parameter $\tau$ denotes the duration between the initial time $\mathring{t}$ and the peak time of a translated shock waveform $T(w(t))$, so the translated peak time is $\mathring{t}+\tau$.

\paragraph{Phase $\varphi$}
This is a trivial parameter since it does not offer any valuable mechanical information, but it is still necessary when a shock signal is decomposed numerically.
This can be considered as a fine tuning of the temporal structure.

\subsection{Connection between the Proposed Shock Waveform and Other Basis Signals}\label{section_waveform_comparision}

The shock waveform proposed in this study can cover a wide range of existing waveforms.
Table \ref{special_cases} shows six special cases of the proposed shock waveform, and their corresponding conditions.
If damping does not exist ($\zeta=0$), the proposed shock waveform is a harmonic wave with frequency $\omega$ and phase $\varphi$.
The Prony mode (damped harmonic wave) is the same as the proposed shock waveform, if damping influence is included and the peak time is the initial time (i.e. $\tau=0$).
Kern and Hayes proposed their waveform by multiplying Prony mode by time variable $t$\cite{kern1984}, which is also a special case of shock waveform when $n=\tau \omega \zeta = 1$.
In far-field condition ($\kappa\gg0$), the shock waveform looks like an asymmetric wavelet (e.g. B-spline wavelet), whose envelop is symmetric but not its phase.
By further restricting the phase parameter of a far-field shock waveform with Eq.(\ref{maximum_condition}), a symmetric wavelet can be obtained, which is very close to the frequently used Morlet wavelet.

The shock waveform is generally similar to the advanced Prony mode proposed by ESA.
A detailed description for the relationship between the proposed shock waveform and the advanced Prony mode is shown in \ref{discussion_AP_SW}.
It shall be noted that, the advance Prony mode denotes those waveforms similar to the far-field shocks in this study, although it can generate a wide range of existing waveforms. 

\begin{table}
	\caption{Special Case of Shock Wave}
	\label{special_cases}
	\centering
	\resizebox{0.6\textwidth}{!}{
		\begin{tabular}{|p{4cm}|l|p{3cm}|}
			\hline 
			Waveforms & Case in Shock Waveform & Shape \\ 
			\hline
			Proposed Waveform & No restriction on parameters & Including any following shape\\
			
			\hline 
			Harmonic Wave $e^{i\omega t}$ & $\zeta=0$ & \parbox{\linewidth}{\includegraphics[width=\linewidth]{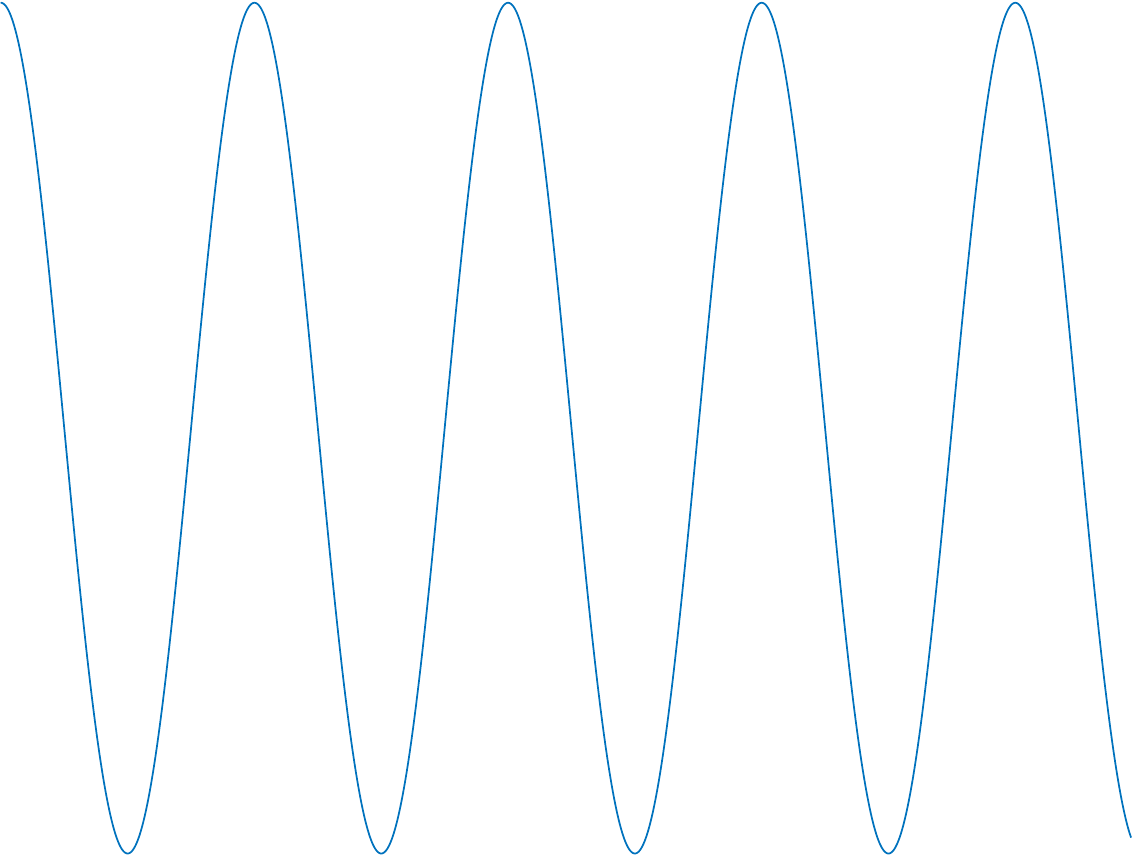}}\\
			
			\hline  
			Prony Mode $e^{i\omega t-\zeta \omega t}$ & $\zeta \neq 0$, $\tau=0$ & \parbox{\linewidth}{\includegraphics[width=\linewidth]{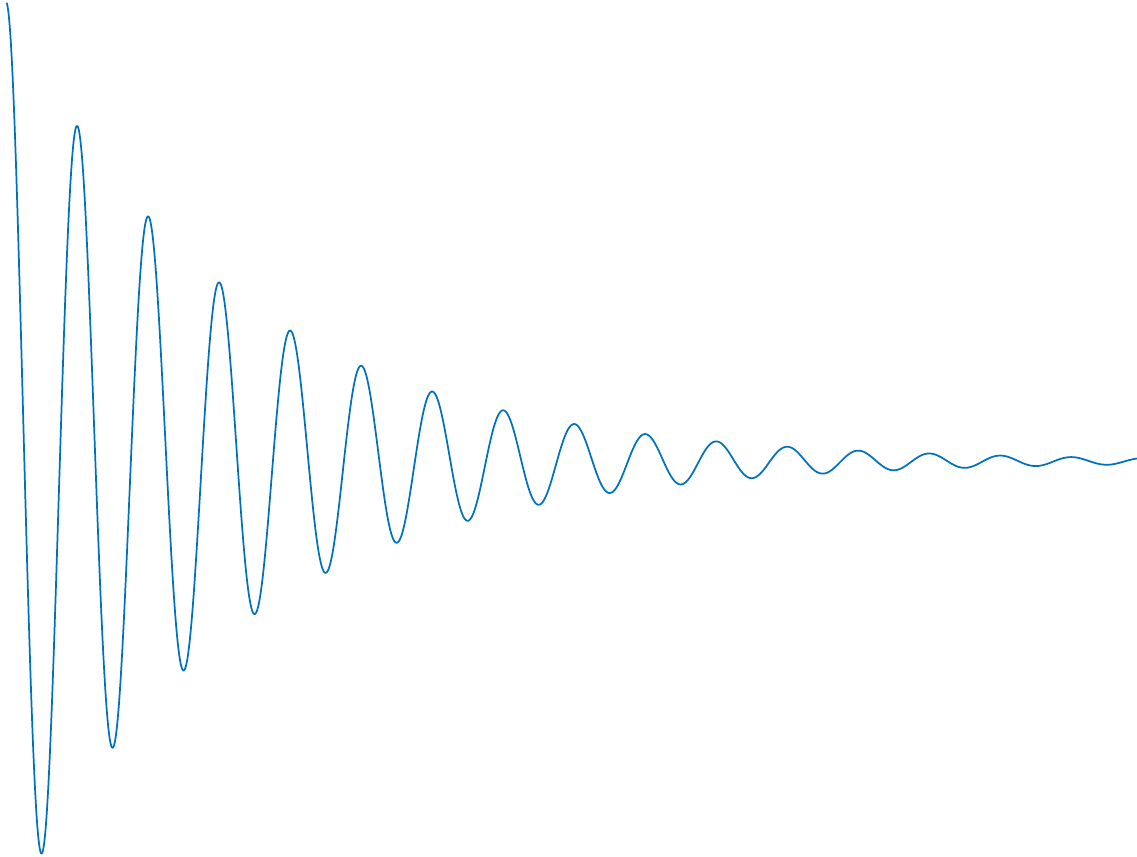}}
			\\
			
			\hline  
			Kern and Hayes' Function $t e^{i\omega t-\zeta \omega t}$ & $n=\tau \omega \zeta = 1$ & \parbox{\linewidth}{\includegraphics[width=\linewidth]{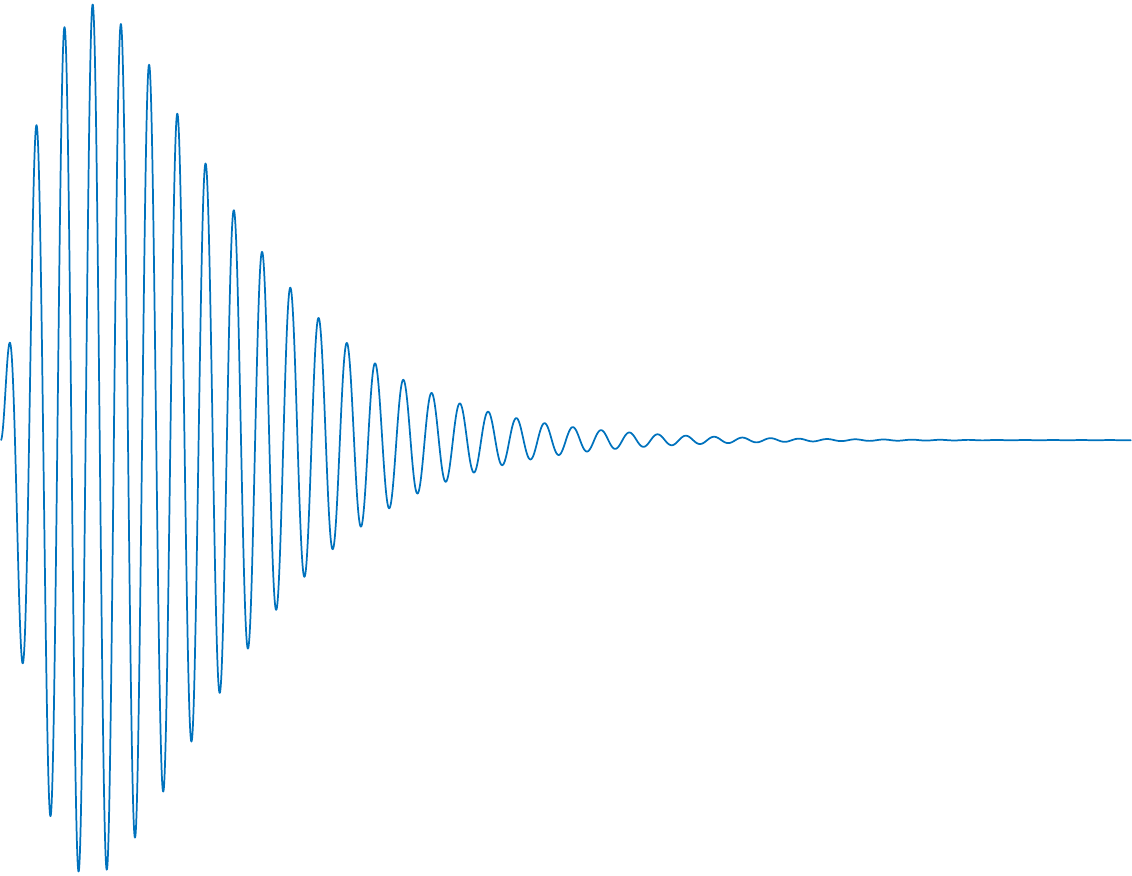}}
			\\
			
			\hline 
			Advanced Prony Mode & $\zeta\neq0, \tau\neq 0$ &
			\parbox{\linewidth}{\includegraphics[width=\linewidth]{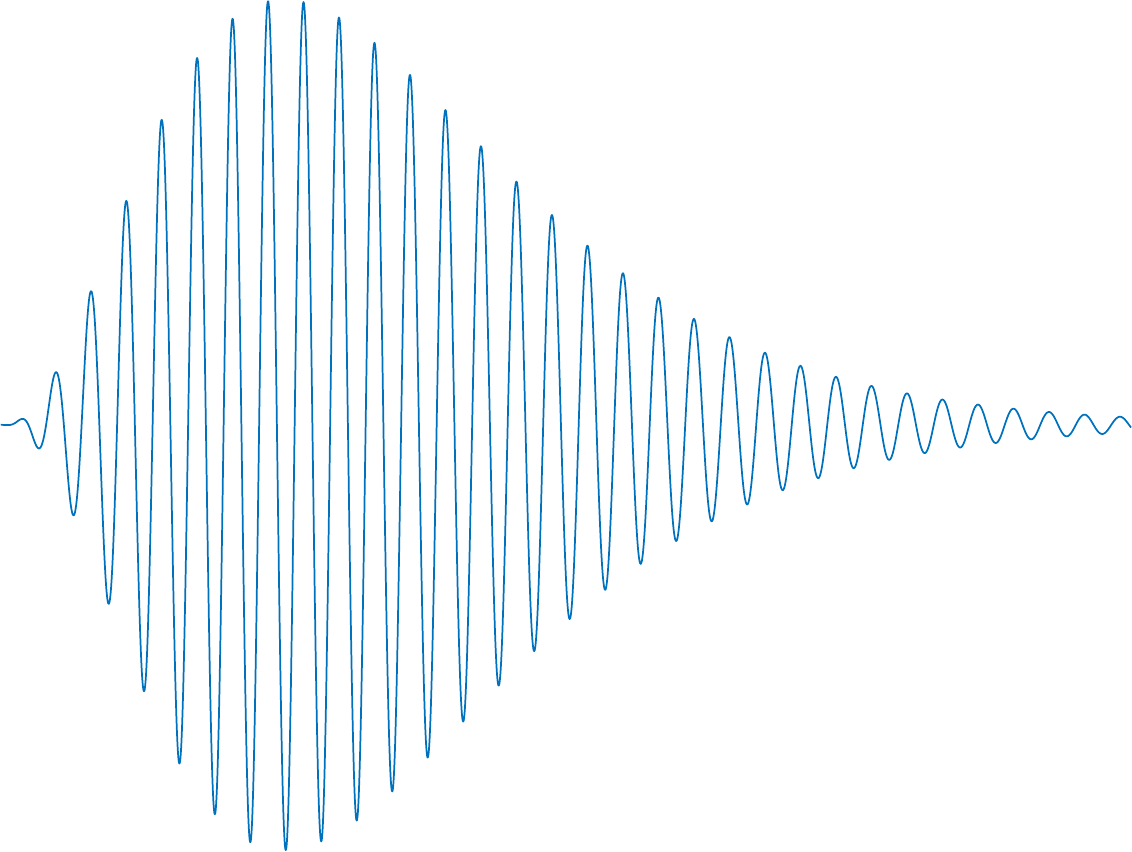}}
			\\ 
			
			\hline 
			Asymmetric Wavelet & $\zeta\neq0, \kappa\gg0$ & \parbox{\linewidth}{\includegraphics[width=\linewidth]{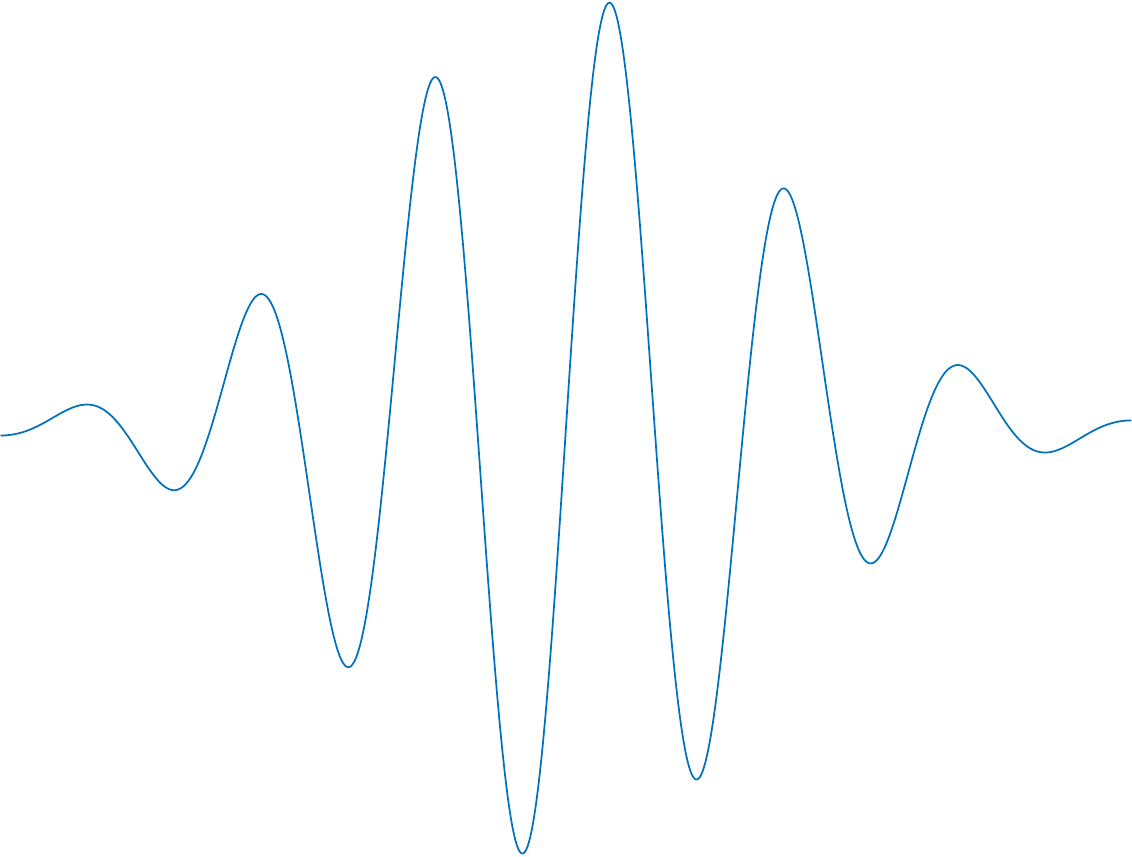}}
			\\
			
			\hline 
			Symmetric Wavelet & $\zeta\neq0, \kappa\gg0, \varphi=-(\tau \omega-\pi \lfloor \frac{\tau \omega}{\pi} \rfloor)$ & \parbox{\linewidth}{\includegraphics[width=\linewidth]{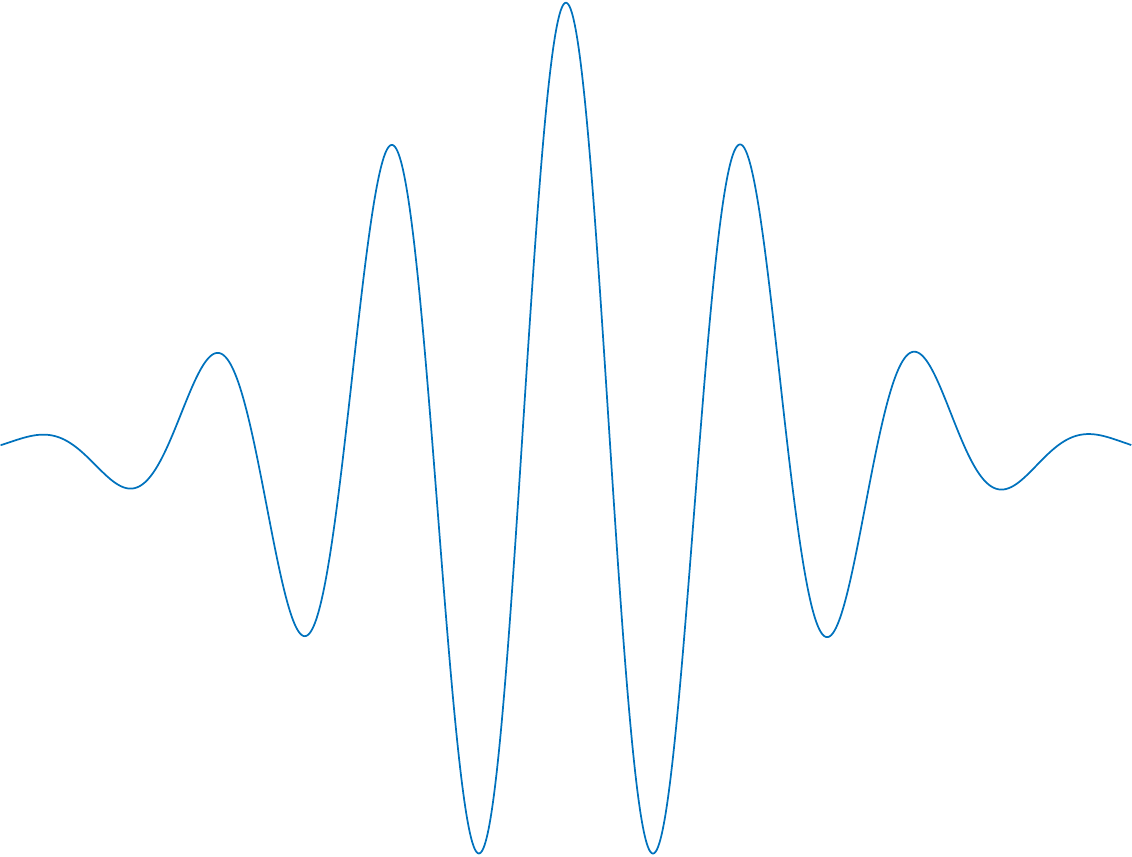}}
			\\
			
			\hline
		\end{tabular} 
	}
\end{table}

\subsection{Fourier Transform of Shock Waveform}

\begin{figure}
	\centering
	\includegraphics[width=0.5\linewidth]{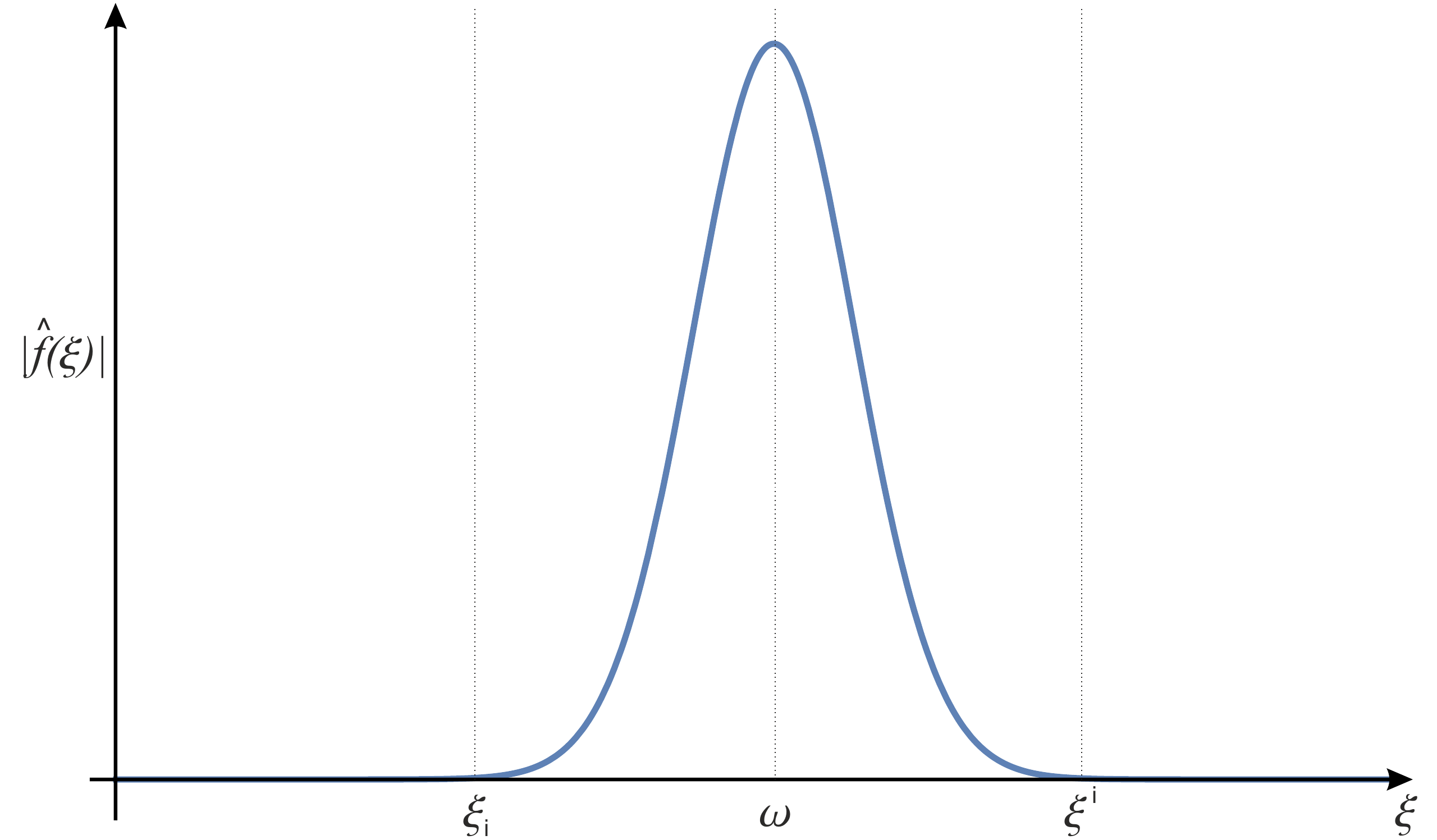}
	\caption{Shape of shock waveform in frequency domain}
	\label{shape_fft}
\end{figure}

The spectrum of the shock waveform in frequency domain can be determined by Fourier transform.
\begin{align}
\hat{f}(\xi) & =\int_{-\infty}^{\infty} w(t) e^{-i \xi t} dt\\
& =\int_{-\infty}^{\infty} At^{\zeta \omega \tau} \tau^{-\zeta \omega \tau} e^{\zeta \omega (\tau-t)+i (\omega t+\varphi)} H(t) e^{-i \xi t} dt\\
& =\int_{0}^{\infty} At^{\zeta \omega \tau} \tau^{-\zeta \omega \tau} e^{\zeta \omega (\tau-t)+i (\omega t+\varphi)} e^{-i \xi t} dt\\
& =A \tau ^{-\zeta  \tau  \omega } (\zeta  \omega +i (\xi -\omega ))^{-\zeta  \tau  \omega -1} \Gamma (\zeta  \tau  \omega +1) e^{\zeta  \tau  \omega +i \varphi }
\end{align}
where $\xi$ represents angular frequency and $\Gamma$ is the Gamma function.
The amplitude of each frequency component is
\begin{equation}\label{fft_amplitude}
|\hat{f}(\xi)|=A e^{\zeta  \tau  \omega } \tau ^{-\zeta  \tau  \omega } (\zeta ^2 \omega ^2+(\xi -\omega )^2)^{\frac{1}{2} (-\zeta  \tau  \omega -1)} \Gamma (\zeta  \tau  \omega +1)
\end{equation}
which satisfies
\begin{equation}
|\hat{f}(\omega+\xi)|=|\hat{f}(\omega-\xi)|
\end{equation}
Hence, the amplitude of frequency component is symmetric about frequency $\omega$, and reaches its peak when $\xi=\omega$.
Shock waveform has a bell-like shape in frequency domain, as given by Eq.(\ref{fft_amplitude}) and depicted in Fig.\ref{shape_fft}.

According to the mechanical mechanism described in \ref{mechanical_background}, this study assumes that each bell spectrum in frequency domain of a shock signal ($r$) can be represented by a shock waveform $w_i$ component in time domain.
If there exists a frequency band $[\xi_i,\xi^i]$ for the $i$th bell, then the following relationship between the band-pass filtered signal, i.e., $r(t)|_{\xi_i}^{\xi^i}$, and a shock waveform is expected.
\begin{flalign}
&&w_i(t)&=r(t)|_{\xi_i}^{\xi^i}& \label{bandpass_relation}\\
\text{with}&&\omega_i&=\frac{\xi_i+\xi^i}{2}&  \label{symmetry_relation}
\end{flalign}
and $\omega_i, A_i, \zeta_i, \mathring{t}_{i}$ and $\tau_i$ are determined by both temporal and frequency structures of $r(t)|_{\xi_i}^{\xi^i}$; $\varphi_i$ can only fine tune the temporal structure in this case.
Examples of Eq.(\ref{bandpass_relation}) will be shown in Section \ref{paragraph_FFS}.

\section{Shock Decomposition Method}

\subsection{Program Implementation}\label{section_program}

\begin{figure}
	\centering
	\includegraphics[height=0.46\textheight]{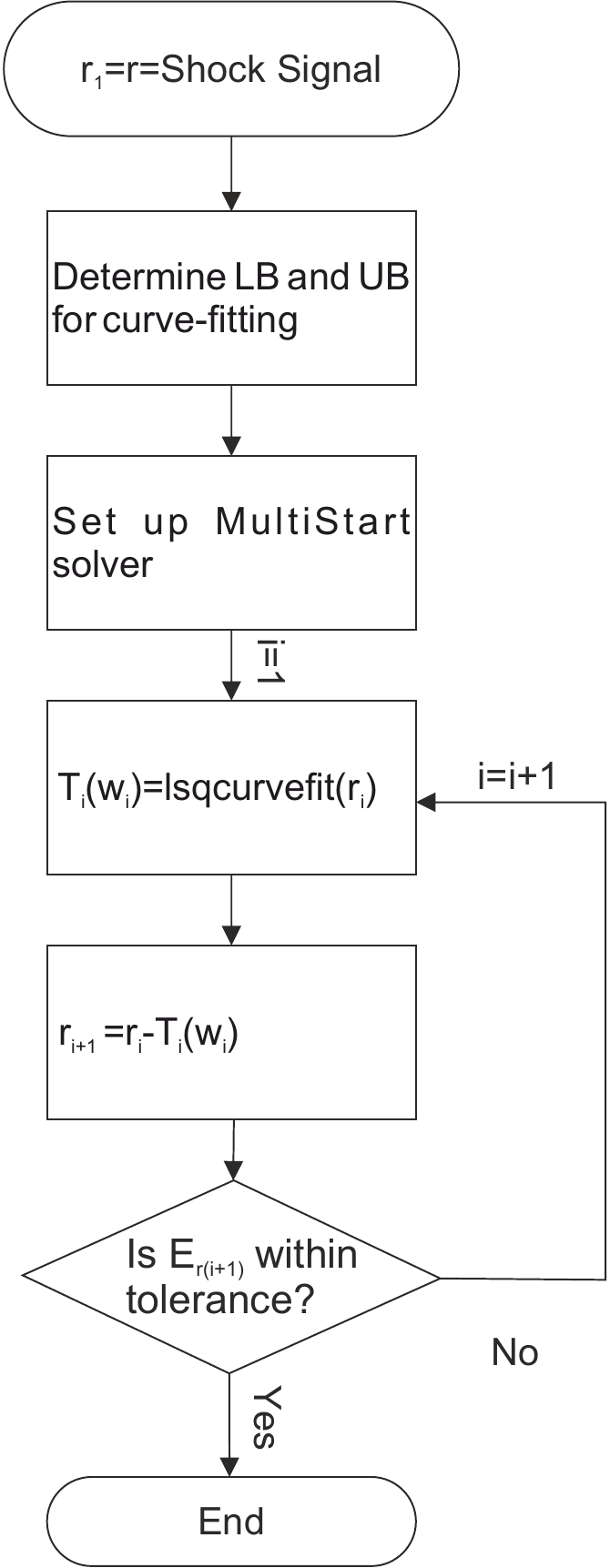}
	\caption{Flowchart of shock decomposition method}
	\label{flowcharts}
\end{figure}

A signal processing tool is proposed to decompose any shock signal into a sum of shock waveform components.
To prevent potential arithmetic overflow during exponentiation calculation, Eq.(\ref{shock_waveform}) is rewritten in the following form,
\begin{equation}\label{waveform_logform}
w(t)=A e^{\zeta\omega\tau(\ln t-\ln \tau )+ \zeta \omega (\tau-t)+i (\omega t+\varphi)} H(t)
\end{equation}
Generally, the main idea of shock decomposition algorithm as shown in Fig.\ref{flowcharts}, is to fit a series of shock waveform components described by Eq.(\ref{waveform_logform}) to a shock signal.
The energy ratio of a waveform component ($w_i$) is calculated in Eq.(\ref{signal_energy}) to evaluate its weight of contribution,
\begin{equation}\label{signal_energy}
\epsilon_{wi}=\frac{E_{wi}}{E_r}
\end{equation}
where $E_{wi}$ is the signal energy of $w_i$ calculated by
\begin{equation}
E_{wi}=\langle w_i(t),w_i(t) \rangle=\int_{-\infty}^{\infty} |w_i(t)|^2 dt
\end{equation}
and $E_{r}$ is the signal energy of shock signal $r$.

The first shock waveform component $T_1(w_1)$ can be obtained by fitting the original shock signal $(r)$, i.e., $r_1=r$. 
The residue, $r_2$, comes from separating the first shock waveform component from the original shock signal by $r_2=r_1-T_1(w_1)$, whose energy ratio $\epsilon_{r2}$ is still not negligible.
This non-linear curve fitting problem
\begin{equation}\label{fitting_problem}
\underset{\x_1}{\text{min}}(E_{r2})=\underset{\x_1}{\text{min}} \int_{-\infty}^{\infty} |r_1(t)-T_1(w_1(t))|^2 dt
\end{equation}
is solved in least-squares sense with `lsqcurvefit' algorithm in Matlab, where $\x_1$ is the solution vector of all parameters to be calculated. 
\begin{equation}
\x_1=[A_1,\omega_1,\mathring{t}_1,\tau_1,\zeta_1,\varphi_1]^{\intercal}
\end{equation}
There may be many local minima of this non-linear problem, in contrast to linear curve fitting where the local minima is also the global minima.
The global solution can be found by using various starting points ($\mathring{\x}_1$) for the fitting, which is achieved by `MultiStart' strategy in Matlab.
The lower bound (LB) and upper bound (UB) of fitting parameters can be set according to the mechanical condition of each shock signal, e.g., peak amplitude, interesting frequency range, etc.
These starting points are discrete combination of parameters in the feasible region defined by LB and UB.

Normally $\epsilon_{w1}$ can account for the main proportion of the shock signal $r$, but the energy ratio $\epsilon_{r2}$ still accounts for non-negligible energy. 
Hence $r_2$ is treated as the new curve to be fitted.
This non-linear fitting procedure can be repeated for all the subsequent residues $r_i$ ($i\geqslant 2$), with $r_{i+1}=r_i-T_i(w_i)$.
\begin{equation}\label{accumulation}
\begin{split}
r2 &= r1-T_1(w_1)\\
r3 &= r2-T_2(w_2)\\
& \enspace\vdots\\
r_{i+1}&=r_i-T_i(w_i)
\end{split}
\end{equation}
\begin{flalign}\label{general_term}
\text{thus,} && r_{i+1}&=r-\sum_{j=1}^{i}T_j(w_j)&
\end{flalign}
According to Eqs.(\ref{accumulation}) and (\ref{general_term}), $r_{i+1}$ is the global error between the shock signal and the reconstructed signal.
The decomposition procedure finally stops when the energy ratio $\epsilon_{r(i+1)}$ of the global error $r_{i+1}$ is within the acceptable tolerance, which is 10\% in this study.

\subsection{Goodness of Fitting and Decomposition}

The goodness of shock waveform decomposition is mainly depended on the set-up of `MultiStart' strategy.
The more starting points, the closer to the global solution (best fit).
But having too many starting points can dramatically increase the computational time.
Therefore, attention should be paid to the choice of feasible starting points.

Basically, starting points shall be decided based on the goodness required.
In this study, starting points are chosen as follow:
\begin{itemize}
	\item $A$: [$max(r)$];
	\item $\omega$: Every octave from $\omega_{LB}$ to $\omega_{UB}$ [$2^0\cdot\omega_{LB},2^1\cdot\omega_{LB},2^2\cdot\omega_{LB},\cdots,\omega_{UB}$];
	\item $\mathring{t}: [-\tau_{r1},\,0,\,\tau_{r1}]$;
	\item $\tau: [\tau_{r1}]$;
	\item $\zeta:[10^{-2},10^{-1},10^{0},10^{1}]$;
	\item $\varphi:[0,\,\frac{\pi}{2},\,\pi,\,\frac{3\pi}{2}]$.
\end{itemize}
Once the starting points are determined, the decomposition process and its results are repeatable. 

The goodness of fitting and decomposition can be evaluated by the energy ratio ($\epsilon_{wi}$) of each shock waveform component $w_i$.
If the solution in each fitting process is global solution, then following relation exists between the values of the energy ratio of consecutive shock waveform components according to Eq.(\ref{fitting_problem}).
\begin{equation}\label{goodness_energy}
\epsilon_{wi}>\epsilon_{w(i+1)}
\end{equation}

The starting points can be set-up based on practical need and engineering experience.
However, it is recommended to increase starting points if $E_{wi}$ check is not satisfied.

\subsection{Selection of Characteristic Shock Waveform Components}

The decomposition process is implemented in least-squares sense, which usually concentrates on the waveform component with higher energy (usually at high frequency).
However, only picking up high energy waveform components may hide low frequency components (usually with less energy), which are also important for the description of shock environment and the analysis of shock response. 
Such a problem can be solved below by the union of two waveform component sets.

Signal energy ratio of the residual signal $\epsilon_{ri}$ is very effective to evaluate the difference of two signals in time domain. 
In this study, the tolerance of energy ratio of the residual signal ($r_{i+1}$) is set to 10\%, which is discarded for its limited influence.
This means that the selected shock waveform components shall account for at least 90\% energy of the shock signal.
The least number of waveform components required to compose 90\% energy of shock signal is defined as $\eta_{90\%}$, which is a parameter to show the complexity of a shock environment.
A set $\{w_i\}_E$ can be used to include these waveform components, i.e.,
\begin{equation}\label{set_energy}
\{w_i\}_E=\{w_i:i=1,2,3,\ldots,\eta_{90\%}\}
\end{equation}

Usually the waveform components set $\{w_i\}_E$ does not cover all important low frequency components.
The method to select important low frequency component is given by
\begin{equation}\label{set_low}
\{w_i\}_L=\{w_i:\omega_i<\omega_j, \; j=1,2,3,\ldots,i-1\},\quad \text{for}\; i\geqslant (\eta_{90\%}+1)
\end{equation}
where all qualified waveform components become the set $\{w_i\}_L$.
Only waveform components with outstanding energy and lower frequency are included in this set.

The united set $\{w_i\}_S$ of $\{w_i\}_E$ and $\{w_i\}_L$ are a set of limited shock waveform components which can represent the time-frequency structure of a shock signal.
\begin{equation}
\{w_i\}_S=\{w_i\}_E \cup \{w_i\}_L
\end{equation}
The simplified shock signal $\hat{r}(t)$ can be reconstructed by waveform components in set $\{w_i\}_S$.
\begin{equation}
\hat{r}(t)=\sum_{i=1} T_i (w_i(t)),\quad \omega_i(t) \in \{w_i\}_S
\end{equation}

\section{Shock Signal Analyses}\label{section_analysis}

Shock waveform can analyse shocks originated from similar mechanical mechanisms, e.g., pyroshock, navy shock and ballistic shock, which are mainly determined by the free vibration of metallic structure (e.g., spacecraft, warship and armoured vehicle).

\begin{figure}
	\centering
	\begin{subfigure}[b]{\textwidth}
		\includegraphics[width=\linewidth]{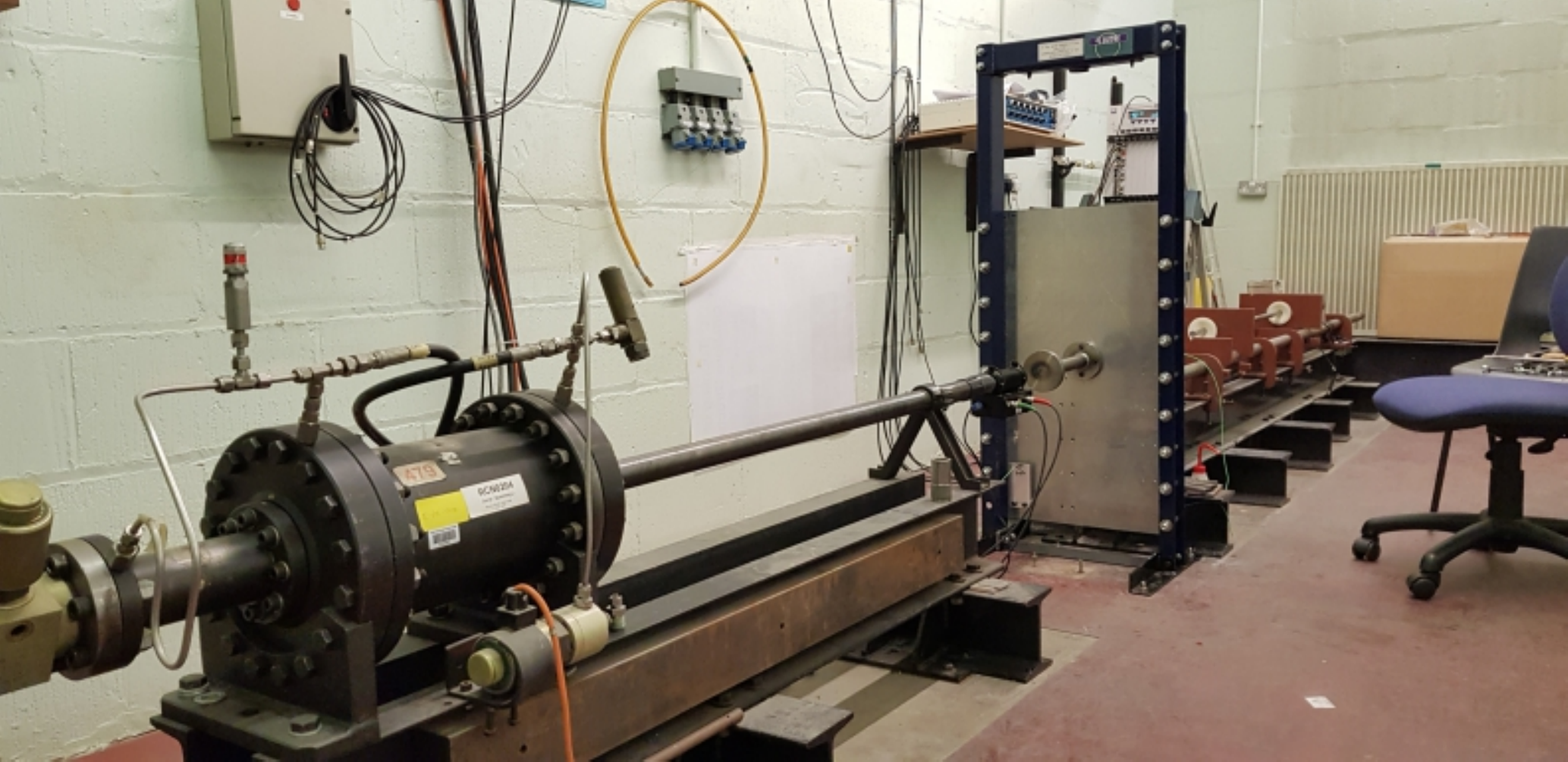}
		\caption{Resonant plate set-up}
		\label{experiment}
	\end{subfigure}
	
	\centering
	\begin{subfigure}[b]{0.49\textwidth}
		\includegraphics[width=\textwidth]{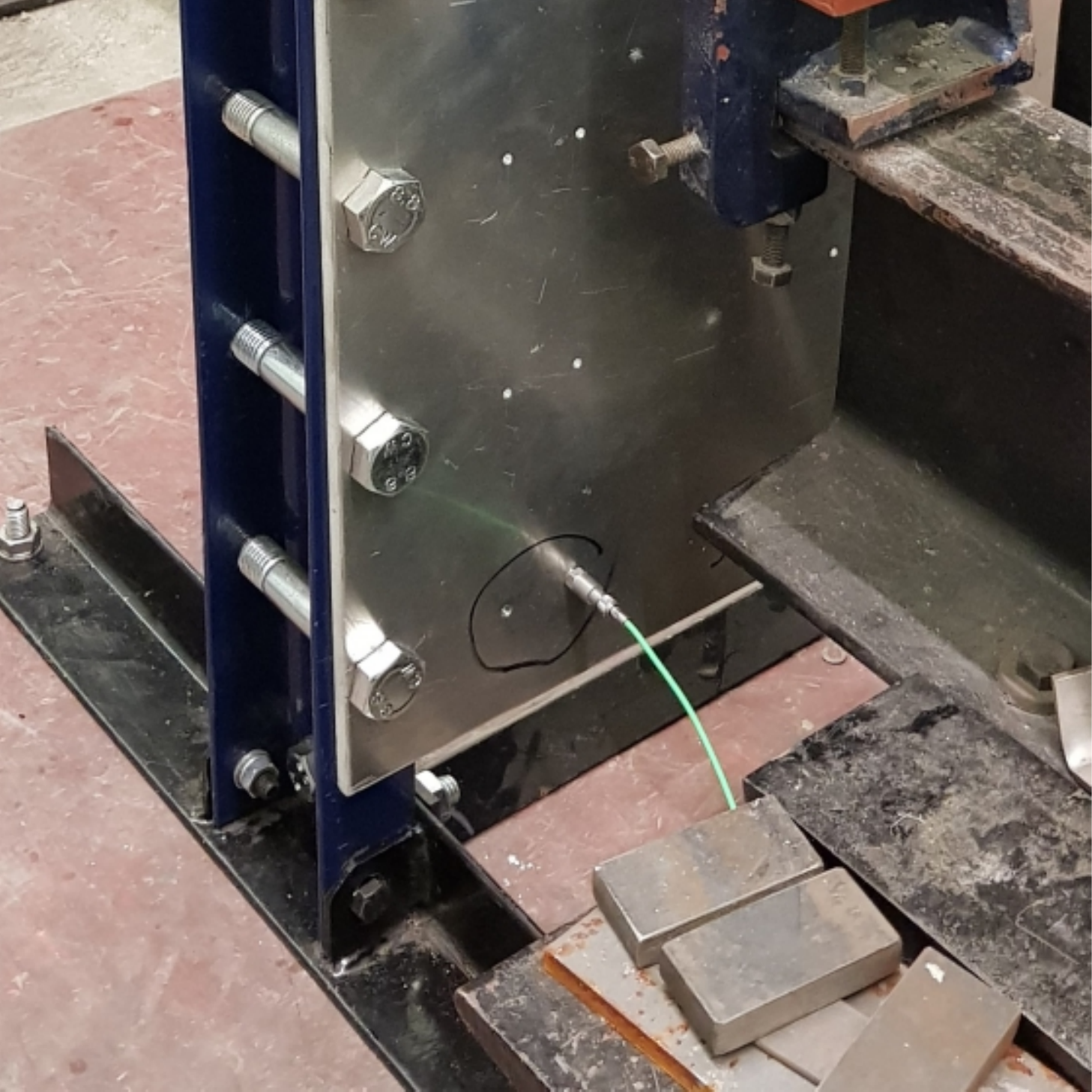}
		\caption{Location of `MFS'}
		\label{location_MFS}
	\end{subfigure}
	\begin{subfigure}[b]{0.49\textwidth}
		\includegraphics[width=\textwidth]{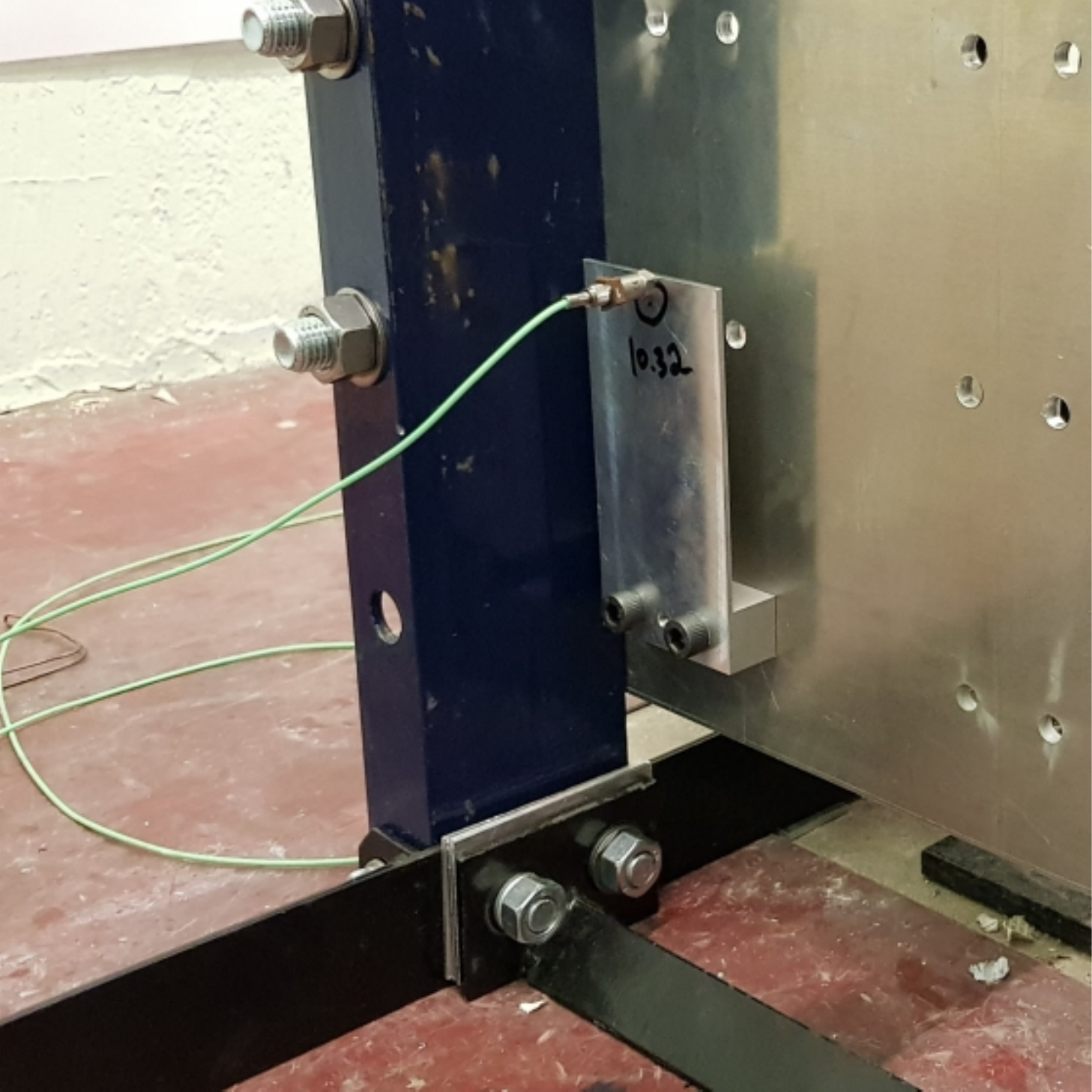}
		\caption{Location of `FFS'}
		\label{location_FFS}
	\end{subfigure}
	\caption{Experimental set-up to produce mid- and far-field shocks simulated by mechanical impact}
	\label{experiment_set_up}
\end{figure}
A near-field pyroshock, measured in the separation stage of an unnamed re-entry rocket vehicle\cite{irvine2010} is taken as an example in this section.
As a comparison, mid- and far-field shocks are produced by mechanical impact (metal-metal impact as shown in Fig. \ref{experiment_set_up} (a)) in laboratory environment where a resonant plate and a sub-structure shown in Fig.\ref{experiment_set_up} (b) and (c), respectively, are used to produce acceleration shock signals for mid- and far-field shocks. 
The shocks are actually produced by the impact of a steel projectile (mass is 58g) and a fully fixed steel resonant plate with a sub-structure, at speed of 13.8m/s.
The accelerometer for data collection is a KISTLER model K-shear 8742A20, whose acceleration and frequency measurement range are within $\pm$20,000g and 100kHz, respectively.

\begin{figure}
	\centering
	\begin{subfigure}[b]{0.32\textwidth}
		\includegraphics[width=\textwidth]{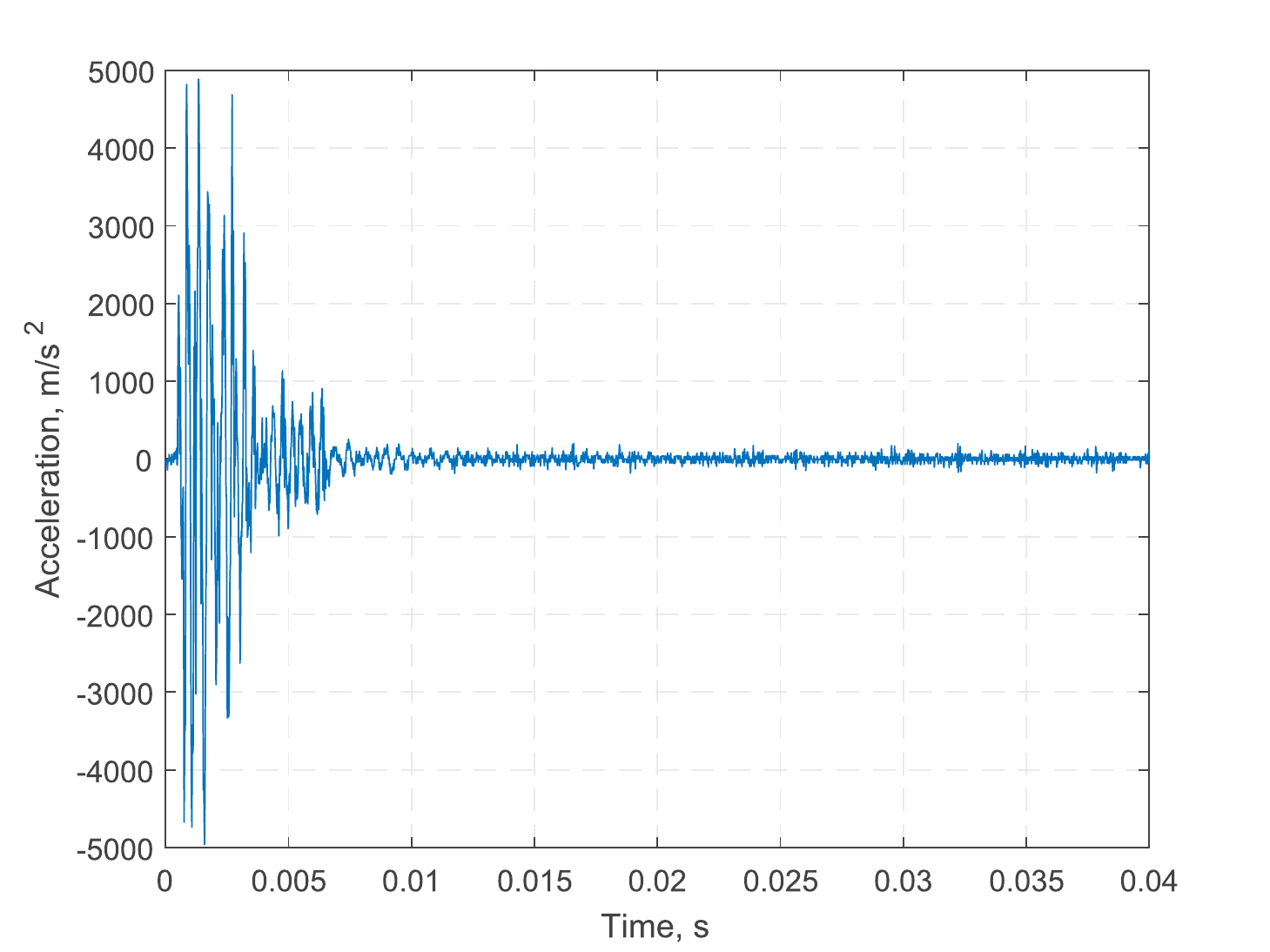}
		\caption{NFS}
		\label{NFS}
	\end{subfigure}
	\begin{subfigure}[b]{0.32\textwidth}
		\includegraphics[width=\textwidth]{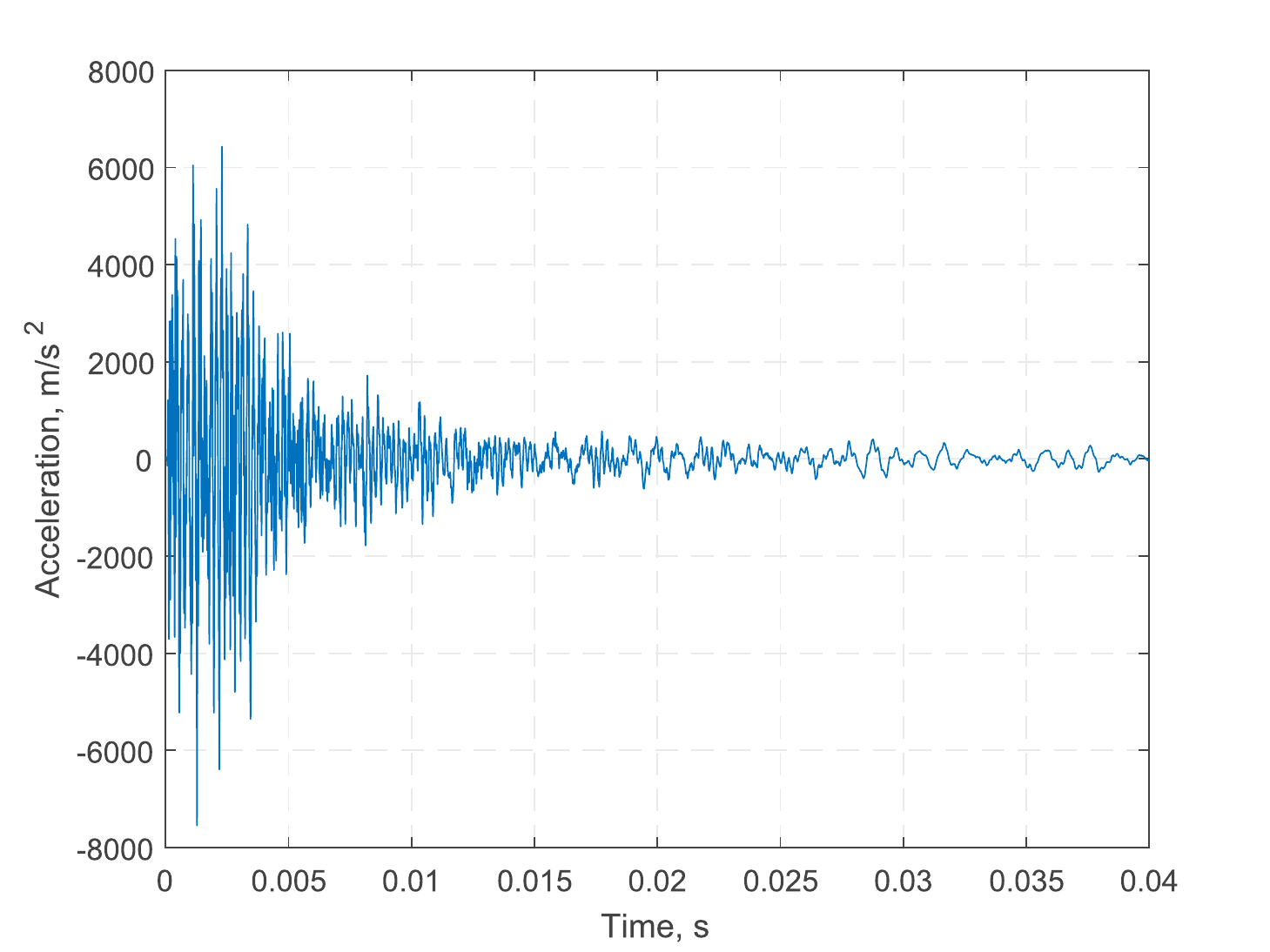}
		\caption{MFS}
		\label{MFS}
	\end{subfigure}
	\begin{subfigure}[b]{0.32\textwidth}
		\includegraphics[width=\textwidth]{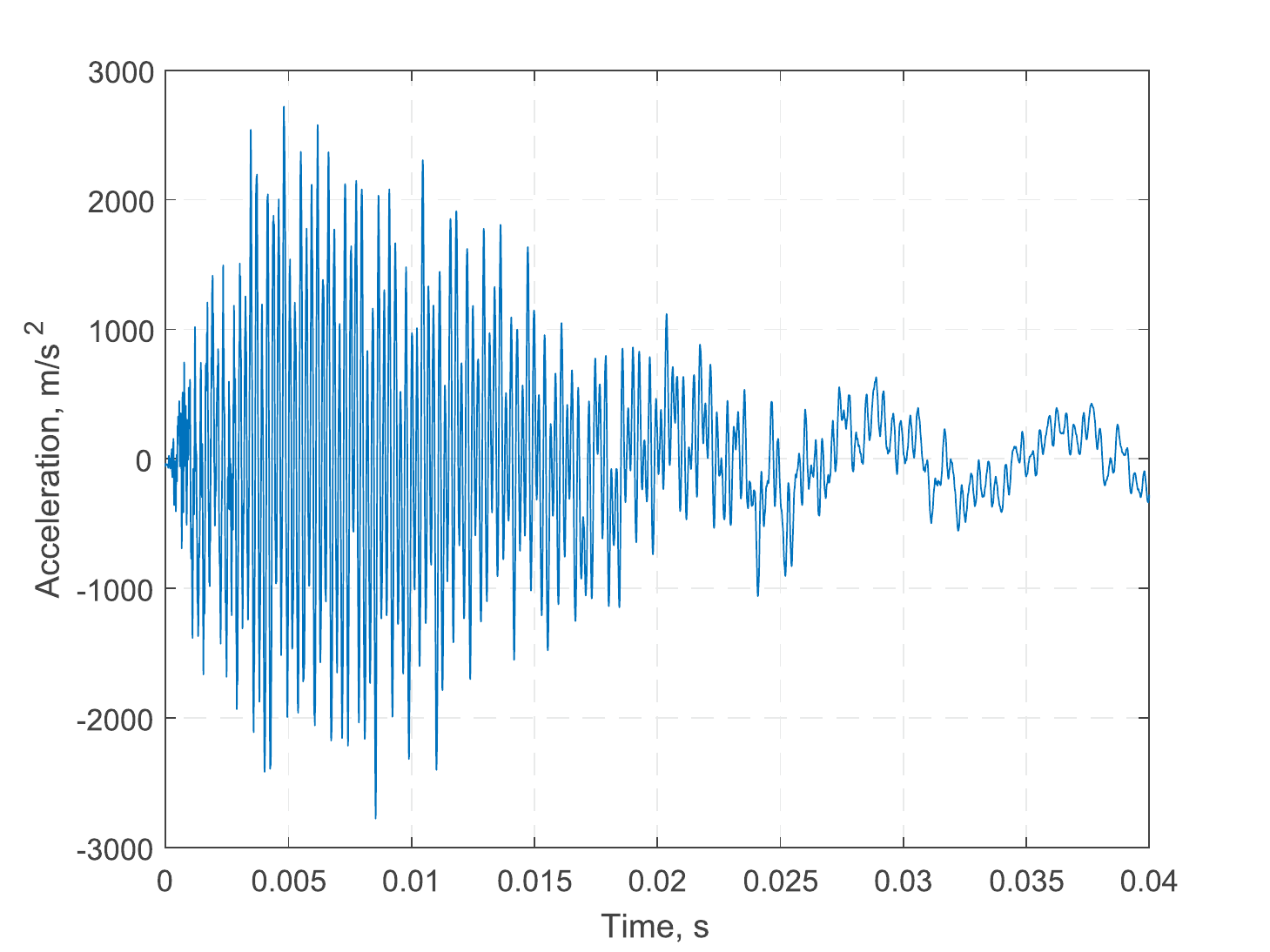}
		\caption{FFS}
		\label{FFS}
	\end{subfigure}
	\caption{Time-histories of NFS, MFS and FFS}
	\label{NFS_MFS_FFS}
\end{figure}

In the subsequent analysis, the near-field pyroshock from re-entry vehicle is termed `NFS'; the mid-field shock from resonant plate and the far-field shock from the sub-structure are referred to `MFS' and `FFS', respectively.
It will be noted that there is no definite distinction between near- and mid-field shocks and between mid- and far-field shocks in common standards.
Here MFS and FFS are considered a `farther' shock than MFS and NFS, respectively, according to the characters described in ECSS Shock Handbook\cite[p.33-35]{ECSS2015}.
All the three shocks are shown in Fig.\ref{NFS_MFS_FFS}.

\subsection{Near-Field Pyroshock from Re-entry Vehicle (NFS)}

\begin{figure}
	\centering
	\includegraphics[width=\linewidth]{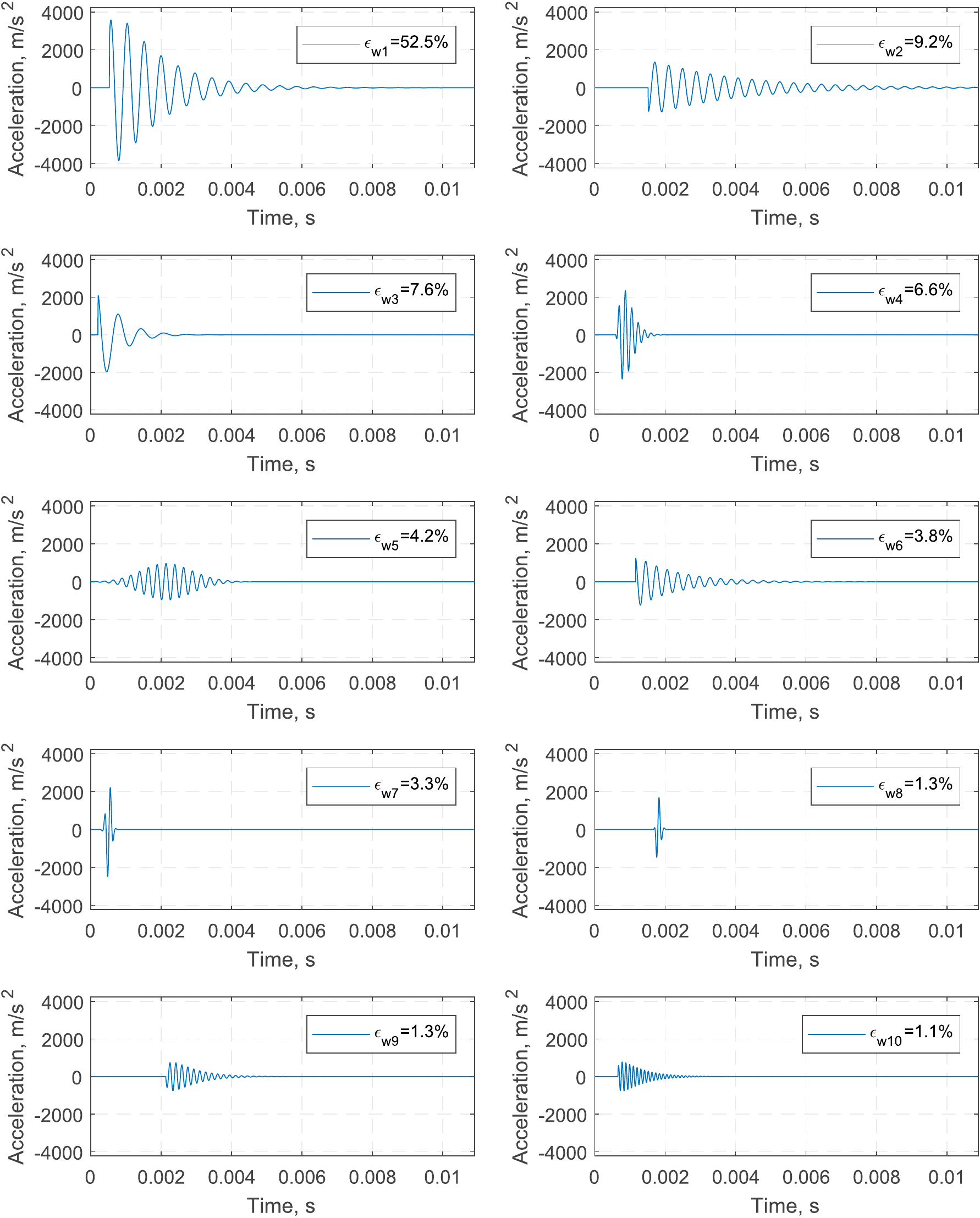}
	\caption{Decomposed shock waveform components of NFS ($T_1(w_1(t))$ to $T_{10}(w_{10}(t))$)}
	\label{w_NFS}
\end{figure}
Fig.\ref{w_NFS} shows the first ten waveform components of NFS from shock waveform decomposition, where each decomposed component of the shock signal is presented in the same scales of acceleration and time.
The $\eta_{90\%}$ of NFS is 10, which means that 90\% of NFS's energy can be represented by the first 10 decomposed components according to the procedure in Section \ref{section_program}.
The goodness of decomposition is acceptable by evaluating the relationship of the signal energy rate of each shock waveform component with Eq.(\ref{goodness_energy}).
The parameters of twelve waveform components in set $\{w_i\}_S$ are listed in Table \ref{xSelect_NFS}.
The first ten components belong to set $\{w_i\}_E$ in energy sense, while the 11th and 12th components are low-frequency components selected by Eq.(\ref{set_low}).

\begin{table}
	\centering
	\caption{Shock waveform components of NFS}
	\label{xSelect_NFS}
	\resizebox{\textwidth}{!}{
		\begin{tabular}{r r r r r r r |r r}
			\hline
			Component&A($m/s^2$)&$\frac{\omega}{2\pi}$(Hz)&$\mathring{t}(ms)$&$\tau(ms)$&$\zeta$&$\varphi$ & $\kappa$ & $\epsilon(\%)$\\
			\hline
			1&3939.41 & 2077 & 0.53 & 0.15 & 0.066 & 5.93 & 0.32 &52.48\\
			\hline
			2&1375.47 & 2534 & 1.52 & 0.08 & 0.027 & 3.31 & 0.21 & 9.22\\
			\hline
			3&2703.28 & 1543 & 0.21 & 0.02 & 0.201 & 0.56 & 0.04 & 7.56\\
			\hline
			4&2435.72 & 5469 & 0.59 & 0.22 & 0.212 & 3.03 & 1.23 & 6.59\\
			\hline
			5&968.94 & 3987 & -10.89 & 13.03 & 1.001 & 0.01 & 51.99 & 4.16\\
			\hline
			6&1367.09 & 3289 & 1.16 & 0.00 & 0.042 & 0.37 & 0.01 & 3.80\\
			\hline
			7&2737.31 & 6056 & -3.15 & 3.66 & 18.987 & 3.18 & 22.19 & 3.28\\
			\hline
			8&1939.25 & 6512 & 1.62 & 0.17 & 1.364 & 3.96 & 1.14 & 1.29\\
			\hline
			9&766.91 & 5896 & 2.13 & 0.18 & 0.050 & 1.89 & 1.07 & 1.27\\
			\hline
			10&781.65 & 9444 & 0.66 & 0.12 & 0.033 & 5.41 & 1.13 & 1.08\\
			\hline
			11&128.05 & 1016 & -10.89 & 12.09 & 1.845 & 3.30 & 12.29 & 0.09\\
			\hline
			12&62.27 & 0.08 & 0.49 & 0.04 & 3333 & 6.22 & 0.00 & 0.01\\
			\hline
		\end{tabular}
	}
\end{table}

The first shock waveform component $w_1$ accounts for more than half of total signal energy, which is obviously a dominant component compared with the rest of waveform components.
The characterisation of temporal structure of the shock signal $r$ can be outlined by $w_1$, whose maximum amplitude $A_1$ is 3939 $m/s^2$ at $\tau_1 + \mathring{t}_1=0.68$ ms.
The frequency $\omega_1$ is 2077 Hz, matching the knee frequency of the SRS or the dominant frequency component of the Fourier spectrum of the shock signal $r(t)$.
Although the damping coefficient of the unnamed re-entry rocket vehicle is not given, the damping ratio $\zeta$ is 6.7\% for $w_1$, close to the general estimation of damping ratio 5\% of spacecraft structural\cite{ECSS2015}.

\begin{figure}
	\centering
	\includegraphics[width=0.5\linewidth]{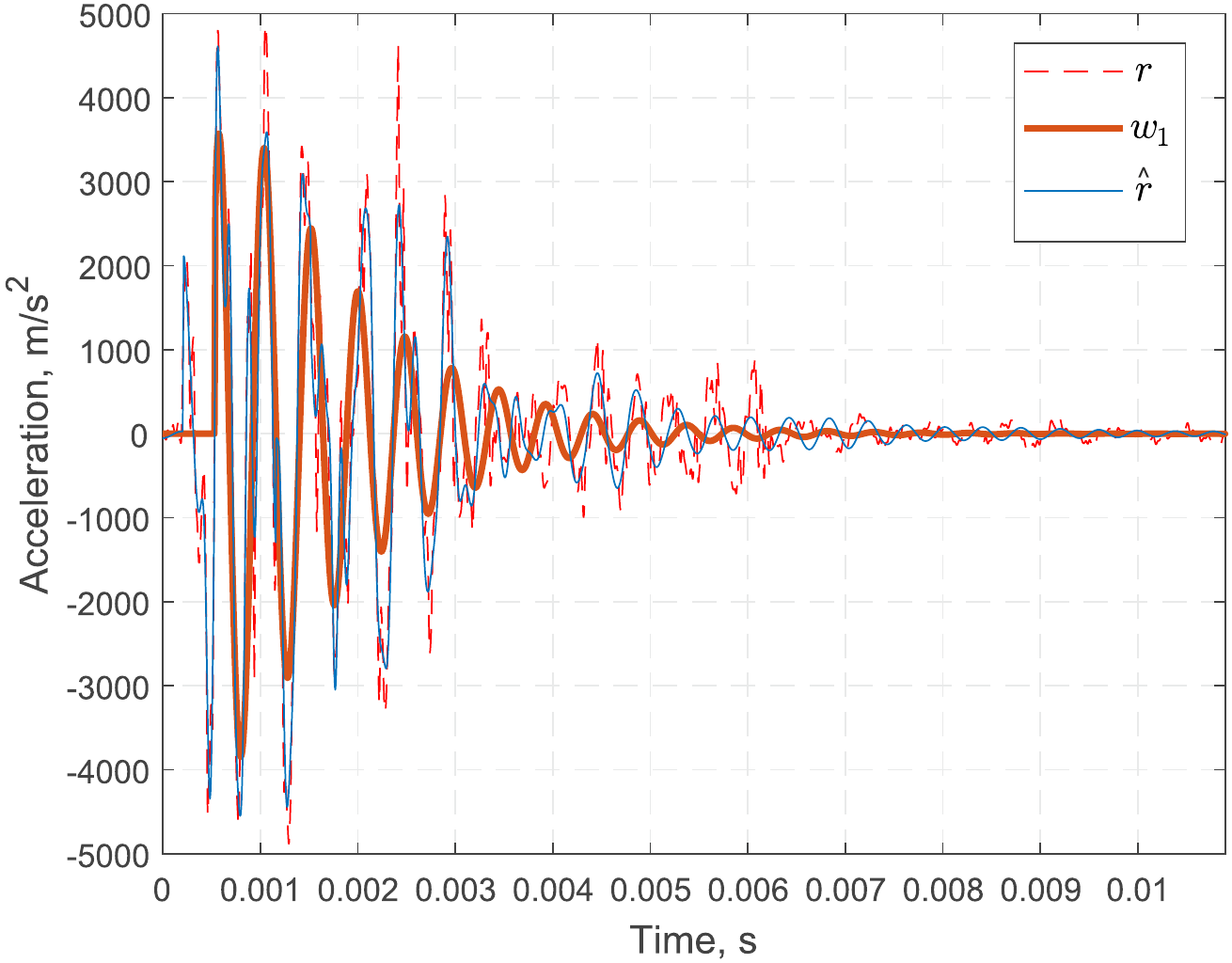}
	\caption{Comparison of $r$, $w_1$ and $\hat{r}$ of NFS}
	\label{r1_w1_sum}
\end{figure}

In terms of the dominant shock distance, the overall shape of this near-field pyroshock ($r$) in Fig.\ref{r1_w1_sum} is visually like the far-field shock rather than the near-field shock shown in Shock Handbook\cite[p35]{ECSS2015}.
Parameter $\kappa$ and $\epsilon$ for each waveform component are given on the right of the vertical line in Table \ref{xSelect_NFS}.
The value $\kappa_1=0.32$ corresponds to a near-field shock scenario and the shape of $w_1$ also belongs to an ordinary damped harmonic wave, which looks like measured near-field shock signal\cite[p.403]{ECSS2015}.
Because $\epsilon_1$ overpowers the rest of $\epsilon_i$, the dominant shock distance can be indicated by $\kappa_1$.

The shock decomposition method can simplify the description of a shock environment while keeping its main temporal structure and frequency content.
The reconstructed waveform $\hat{r}=\sum_i \{w_i\}_S$ shall still represent the shock signal $r$ in both time and frequency domain.
The time-history comparison of $r$, $w_1$ and reconstructed waveform $\hat{r}$ is illustrated in Fig.\ref{r1_w1_sum}.
The 10\% energy ratio of the residual signal $(r-\hat{r})$ guarantees the similarity between $r$ and $\hat{r}$ in time domain.

\begin{figure}
	\centering
	\begin{subfigure}[b]{0.45\textwidth}
		\includegraphics[width=\textwidth]{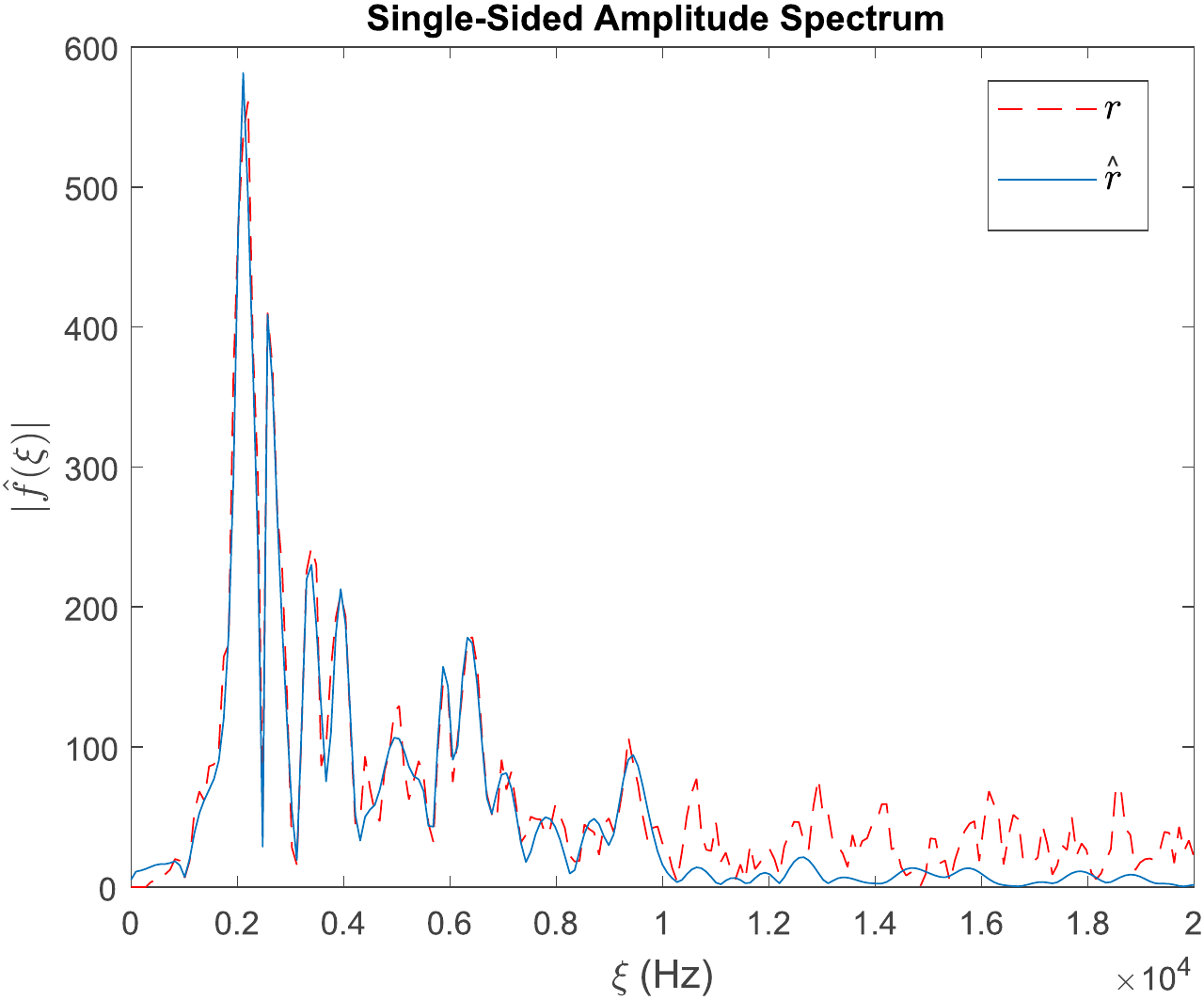}
		\caption{FFT spectrum}
		\label{fft_NFS}
	\end{subfigure}
	\begin{subfigure}[b]{0.45\textwidth}
		\includegraphics[width=\textwidth]{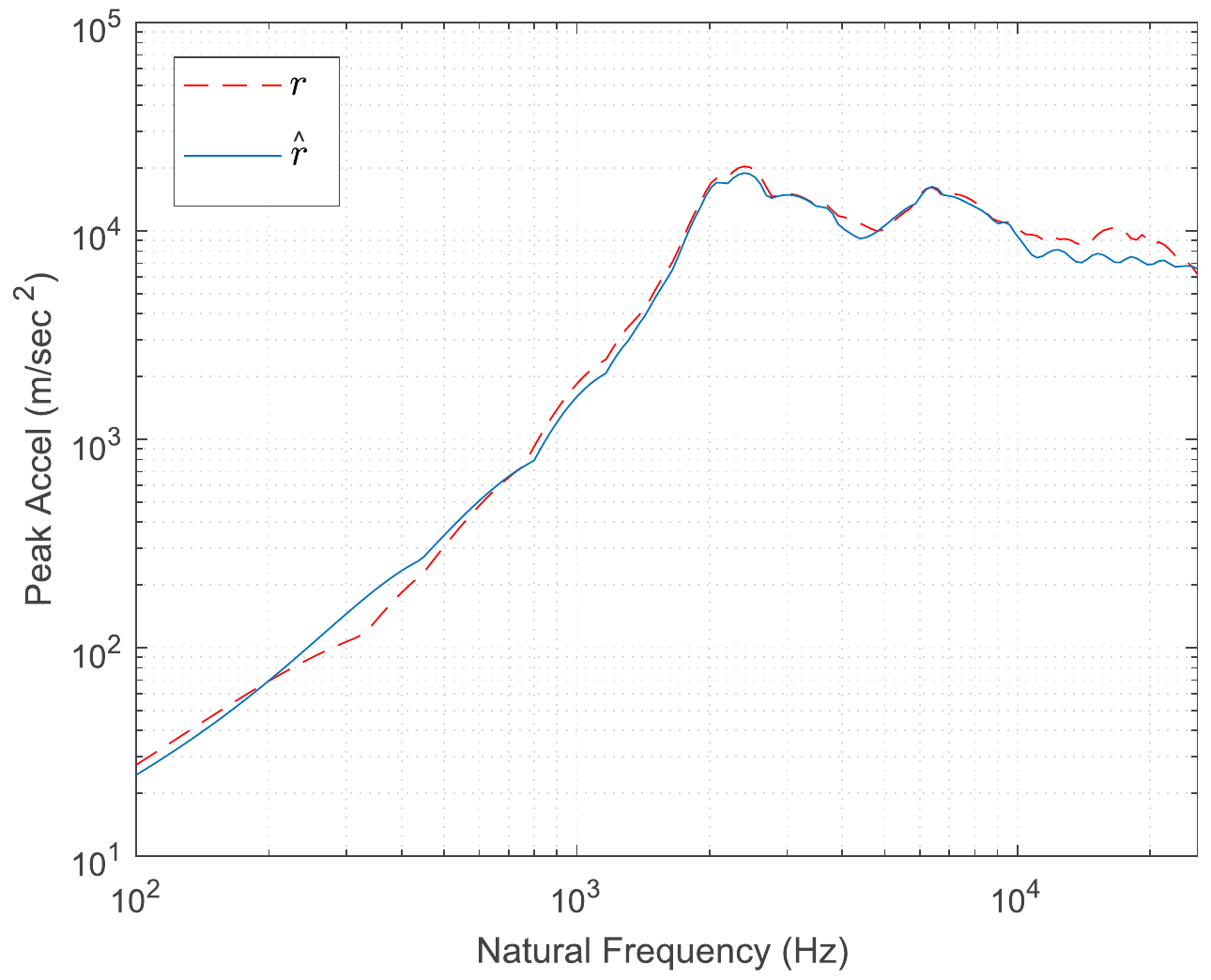}
		\caption{SRS}
		\label{srs_NFS}
	\end{subfigure}
	\caption{FFT and SRS of $r$ and $\hat{r}$ of NFS}
	\label{fft_srs_NFS}
\end{figure}
Comparisons of $\hat{r}$ and $r$ in frequency domain are presented in two forms, i.e., frequency spectrum in Fig.\ref{fft_srs_NFS}(a) and SRS in Fig.\ref{fft_srs_NFS}(b).
The spectrum of $\hat{r}$ can cover several most significant bells of the spectrum of $r$ in frequency domain, depending on how many shock waveform components are involved.
In this case, $\eta_{90\%}$ of NFS is ten, which means that ten bells are covered.
The frequencies where maximum errors take place are the centre frequencies of components $w_i$ outside set $\{w_i\}_E$, whose energy ratios are less than 1\%.
These components are discarded since the cumulative error in the whole frequency region is less than the acceptable tolerance, i.e., 10\%.
The SRS of $\hat{r}$ also matches the SRS of $r$ very well.
There exists certain error in low-frequency domain if shock waveform components in set $\{w_i\}_L$ are not included, which shows the importance of low-frequency waveform components for mechanical response.

\subsection{Mid-Field Shock from Resonant Plate (MFS)}

\begin{figure}
	\centering
	\includegraphics[width=\linewidth]{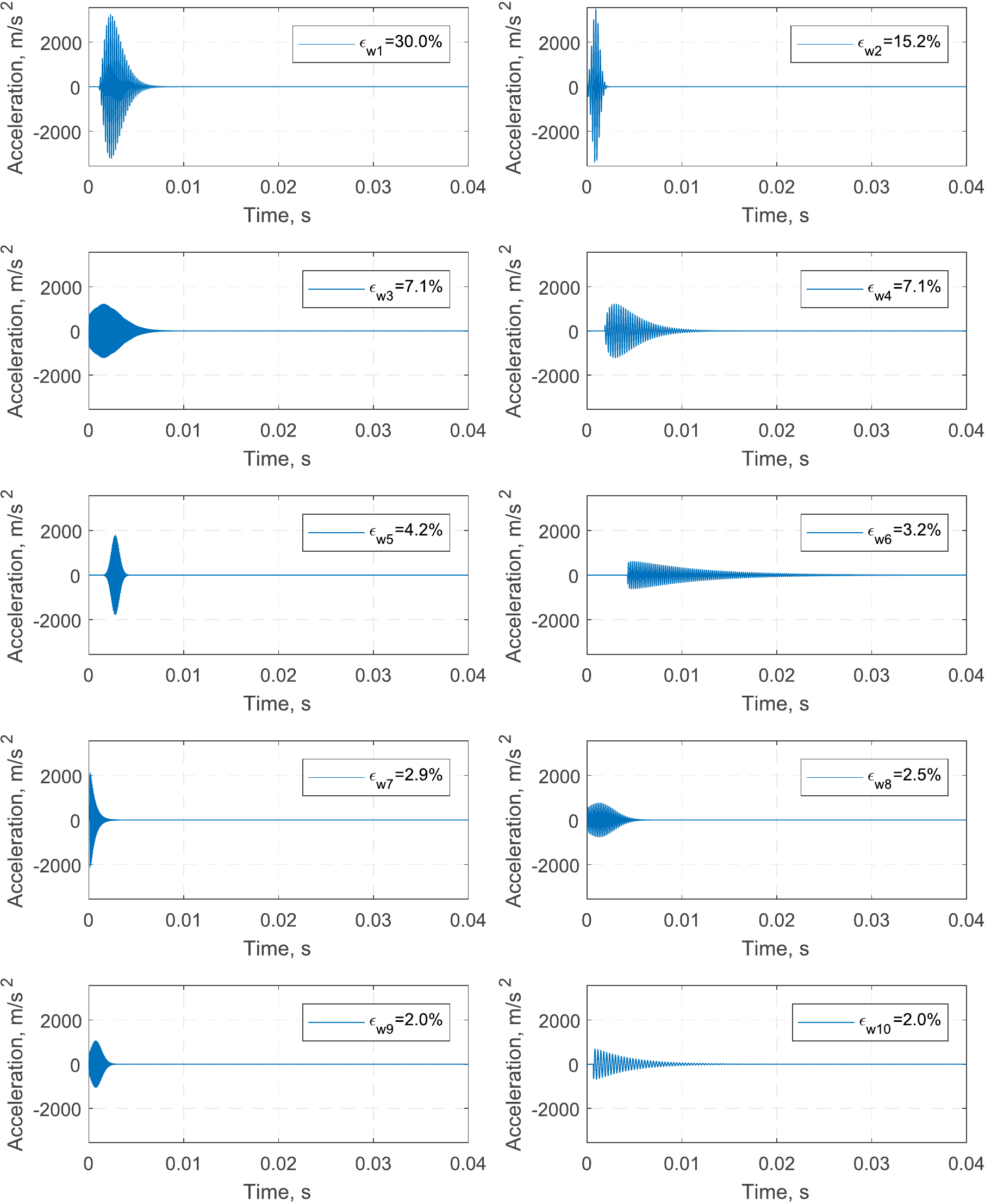}
	\caption{Decomposed shock waveform components of MFS}
	\label{w_MFS}
\end{figure}

The first ten shock waveform components are shown in Fig.\ref{w_MFS}.
A more complicated mechanical scenario is indicated by $\eta_{90\%}$ of MFS having 25 and 3 components in $\{w_i\}_E$ and $\{w_i\}_L$, respectively.
Although 28 waveform components are required in total to characterize this complex shock environment, only parameters for the first ten components in $\{w_i\}_E$ are listed in Table \ref{xSelect_MFS} for illustration purpose.

\begin{table}
	\centering
	\caption{Shock waveform components of MFS}
	\label{xSelect_MFS}
	\resizebox{\textwidth}{!}{
		\begin{tabular}{r r r r r r r | r r}
			\hline
			Component&A($m/s^2$)&$\frac{\omega}{2\pi}$(Hz)&$\mathring{t}(ms)$&$\tau(ms)$&$\zeta$&$\varphi$ & $\kappa$ & $\epsilon(\%)$ \\
			\hline
			1&3231.82 & 4683 & 0.99 & 1.3 & 0.056 & 6.01 & 6.11 & 30.02 \\
			\hline
			2&3489.69 & 4315 & -21.08 & 21.99 & 5.272 & 0.11 & 94.89 & 15.20 \\
			\hline
			3&1206.31 & 22068 & -3.93 & 5.46 & 0.015 & 1.90 & 120.69 & 7.14 \\
			\hline
			4&1222.74 & 4086 & 1.84 & 1.07 & 0.021 & 5.87 & 4.39 & 7.13 \\
			\hline
			5&1786.27 & 16446 & -3.88 & 6.66 & 0.385 & 1.59 & 109.63 & 4.15 \\
			\hline
			6&615.19 & 4863 & 4.25 & 0.53 & 0.005 & 3.44 & 2.60 & 3.21 \\
			\hline
			7&2138.04 & 16626 & 0.09 & 0.04 & 0.021 & 5.46 & 0.71 & 2.89 \\
			\hline
			8&762.57 & 5554 & -18.25 & 19.51 & 0.239 & 6.11 & 108.40 & 2.52 \\
			\hline
			9&1059.31 & 20996 & -3.64 & 4.37 & 0.099 & 0.04 & 91.91 & 2.01 \\
			\hline
			10&-693.35 & 3249 & 0.65 & 0.18 & 0.014 & 5.75 & 0.59 & 1.99 \\
			\hline
		\end{tabular}
	}
\end{table}

In the case of MFS, even though the first shock waveform component is not as outstanding as the first component of NFS, it can still provide some important features of MFS.
The parameter $\kappa_1=6.12$ identifies MFS as a `farther' shock compared to NFS, although they have similar overall shape visually.
This coincides with the description of resonant plate technology\cite{Lee2012}, and hence, the MFS signal is regarded as a mid-field shock in this study.
The damping ratio of the resonant plate apparatus is indicated by its $\zeta_1$ (5.7\%), which coincides with the 5\% general estimation of damping ratio as well.

\begin{figure}
	\centering
	\includegraphics[width=0.5\linewidth]{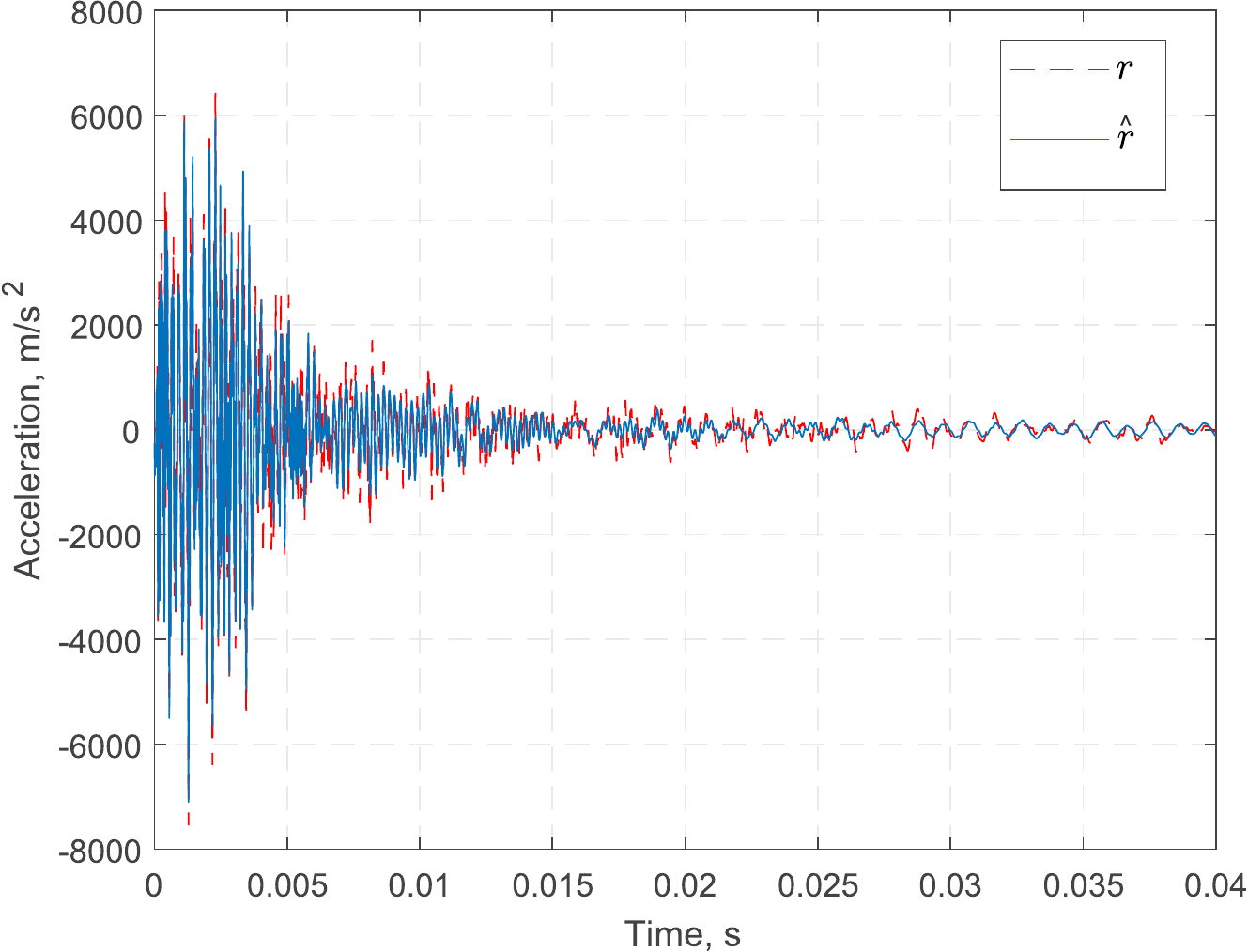}
	\caption{Comparison of $r$ and $\hat{r}$ of MFS}
	\label{acc_r1_sum}
\end{figure}

\begin{figure}
	\centering
	\begin{subfigure}[b]{0.45\textwidth}
		\includegraphics[width=\textwidth]{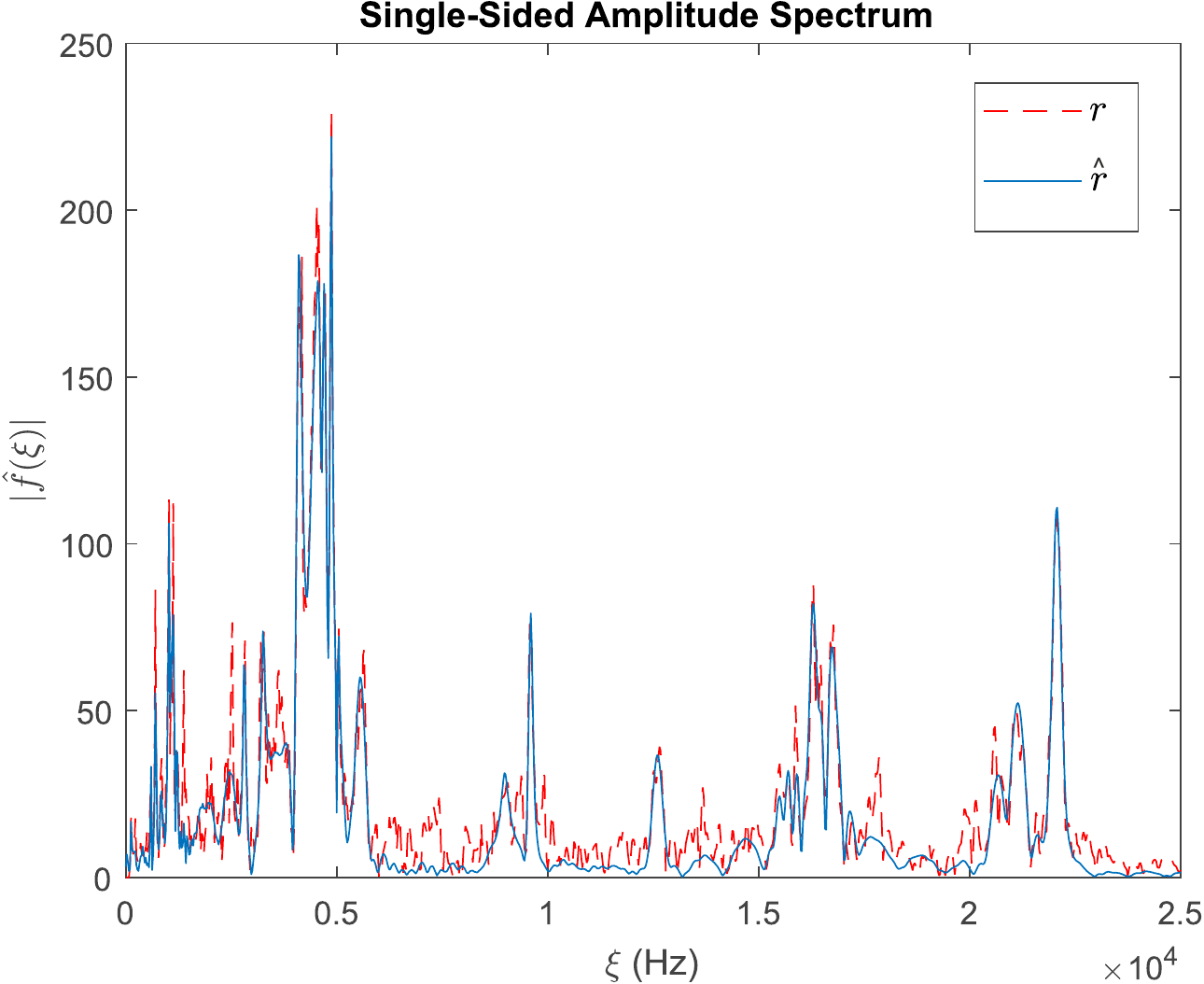}
		\caption{FFT spectrum}
		\label{fft_MFS}
	\end{subfigure}
	\begin{subfigure}[b]{0.45\textwidth}
		\includegraphics[width=\textwidth]{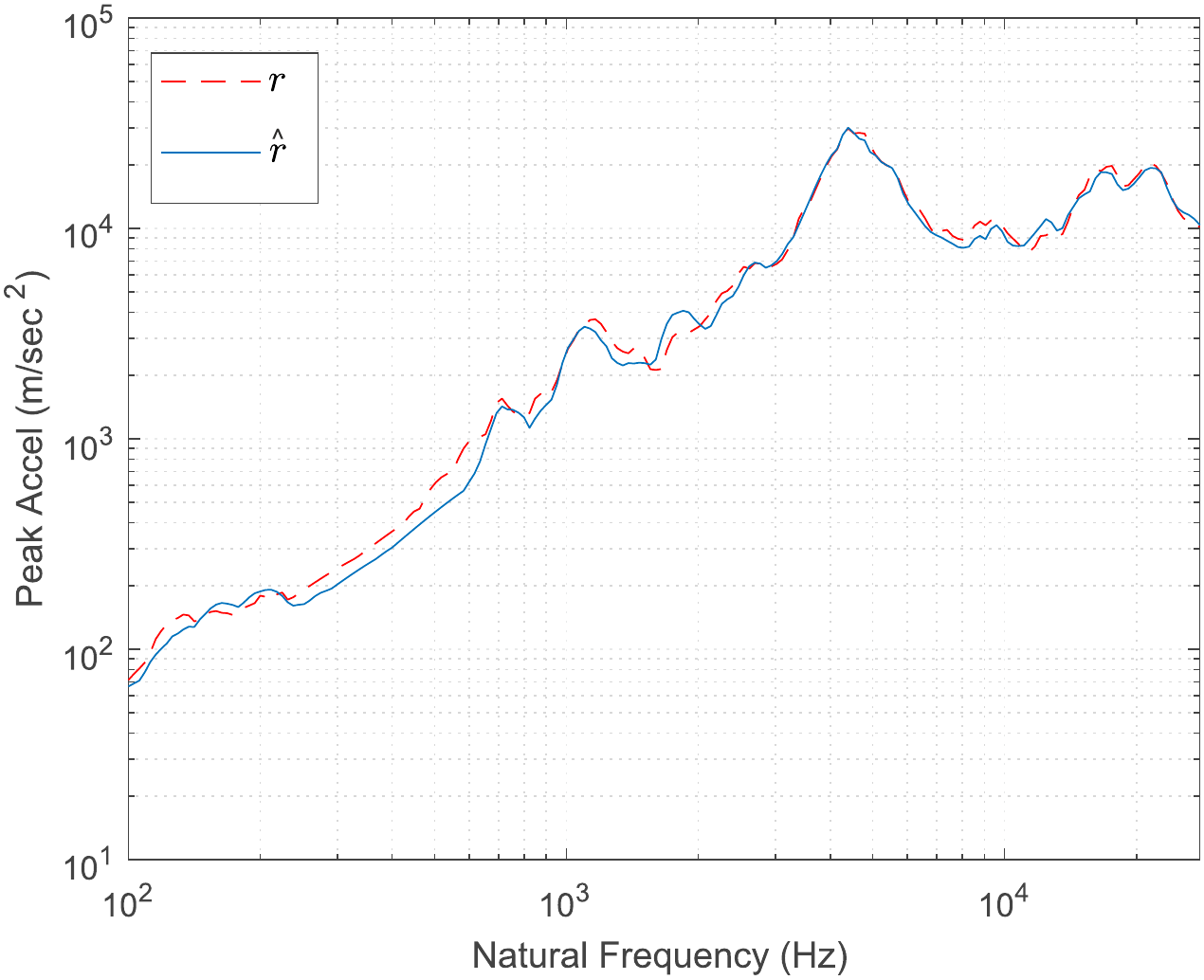}
		\caption{SRS}
		\label{srs_MFS}
	\end{subfigure}
	\caption{FFT and SRS of $r$ and $\hat{r}$ of MFS}
	\label{fft_srs_MFS}
\end{figure}

The comparison of $r$ and $\hat{r}$ in both time and frequency domain are shown in Figs.\ref{acc_r1_sum} and \ref{fft_srs_MFS}, respectively.
The frequency spectrum may also explain the complexity of frequency component of MFS. 

\subsection{Far-Field Shock from Sub-Structure (FFS)}\label{paragraph_FFS}

\begin{figure}
	\centering
	\includegraphics[width=\linewidth]{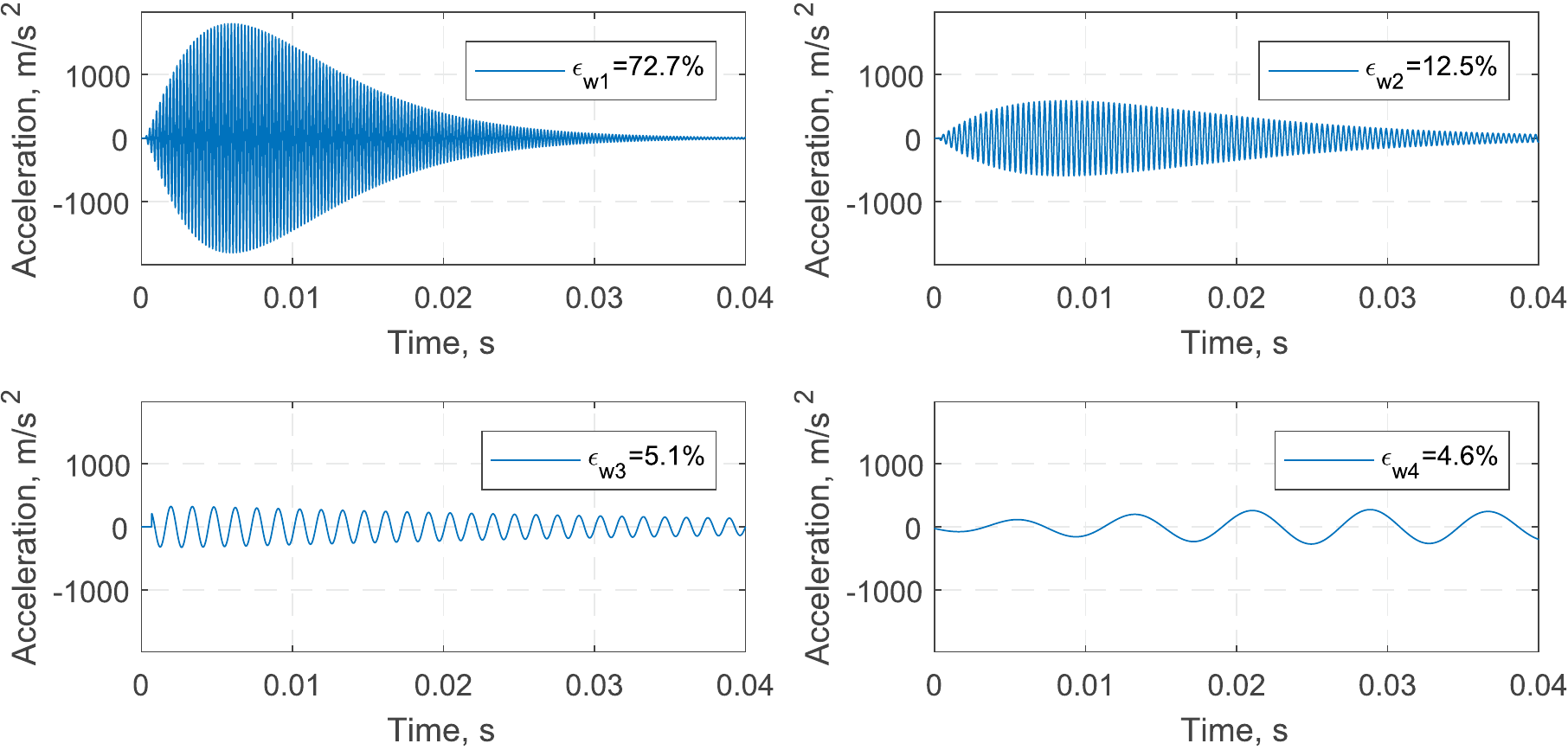}
	\caption{Decomposed shock waveform components of FFS}
	\label{w_FFS}
\end{figure}

This far-field shock scenario is a strong support of the proposed shock waveform as a characteristic waveform.
With only three waveform components, 90\% signal energy of the shock signal can be represented.
Time-history and parameters of $\{w_i\}_S$ are shown in Fig.\ref{w_FFS} and Table \ref{xSelect_FFS}, respectively, where waveform components 4 and 5 belong to $\{w_i\}_L$.

\begin{table}
	\centering
	\caption{Shock waveform components of FFS}
	\label{xSelect_FFS}
	\resizebox{\textwidth}{!}{
		\begin{tabular}{r r r r r r r | r r}
			\hline
			Component&A($m/s^2$)&$\frac{\omega}{2\pi}$(Hz)&$\mathring{t}(ms)$&$\tau(ms)$&$\zeta$&$\varphi$& $\kappa$ & $\epsilon(\%)$ \\
			\hline
			1&1804.18 & 4437 & 0.22 & 5.73 & 0.0079 & 3.63 & 25.45 & 72.73 \\
			\hline
			2&592.20 & 2829 & 0.20 & 8.36 & 0.0071 & 6.11 & 23.66 & 12.53 \\
			\hline
			3&324.35 & 702 & 0.67 & 1.57 & 0.0061 & 0.64 & 1.10 & 5.14 \\
			\hline
			4&274.58 & 127 & -39.99 & 67.21 & 0.2375 & 1.36 & 8.57 & 4.61 \\
			\hline
			5&-32.78 & 90 & 1.44 & 1.00 & 0.0307 & 6.22 & 0.09 & 0.06 \\
			\hline
		\end{tabular}
	}
\end{table}

The first shock waveform component $w_1$ of FFS is particularly important.
It constitutes 72.7\% of total energy of the shock signal.
This shock environment is a common far-field shock, which can be validated both from its temporal shape and the parameter $\kappa_1=25.46$.
The damping ratio of the sub-structure is implied by $\zeta_1=0.79\%$, which is much smaller than the damping of the resonant plate.
This difference may be mainly resulted from different mounting and boundary condition.

\begin{figure}
	\centering
	\includegraphics[width=0.5\linewidth]{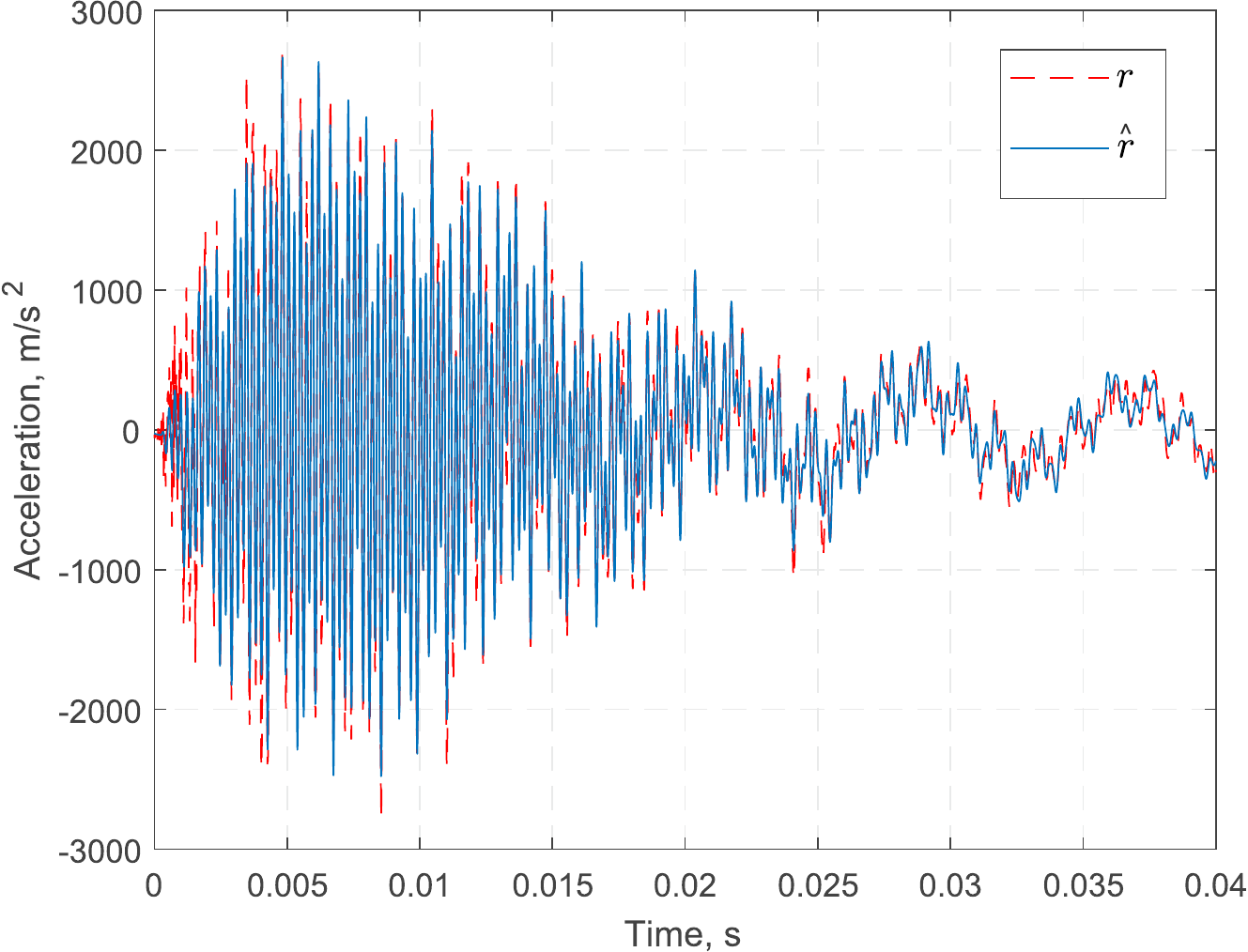}
	\caption{Comparison of $r$ and $\hat{r}$ of \textsc{FFS}}
	\label{acc_FFS}
\end{figure}

\begin{figure}
	\centering
	\begin{subfigure}[b]{0.45\textwidth}
		\includegraphics[width=\textwidth]{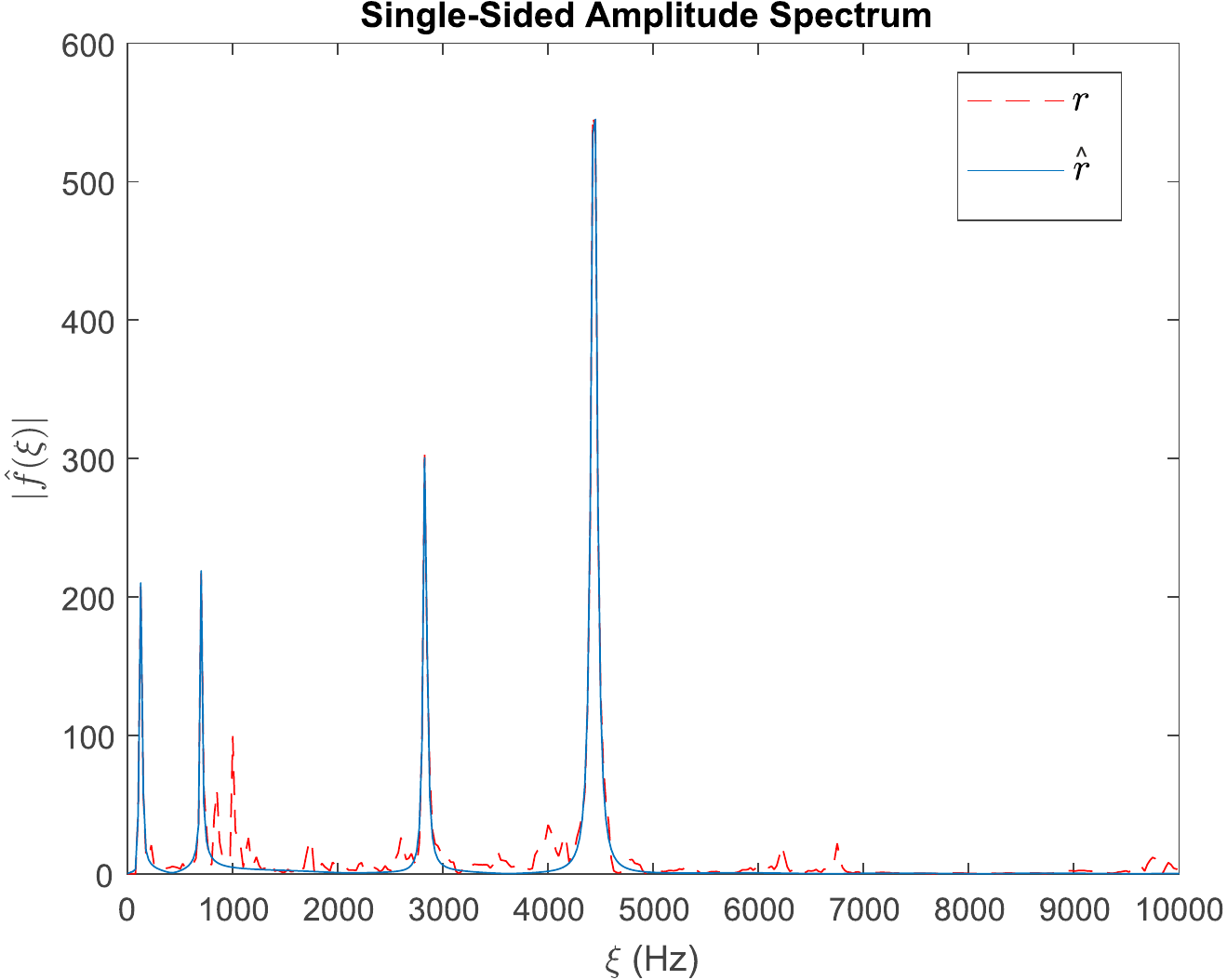}
		\caption{FFT spectrum}
		\label{fft_FFS}
	\end{subfigure}
	\begin{subfigure}[b]{0.45\textwidth}
		\includegraphics[width=\textwidth]{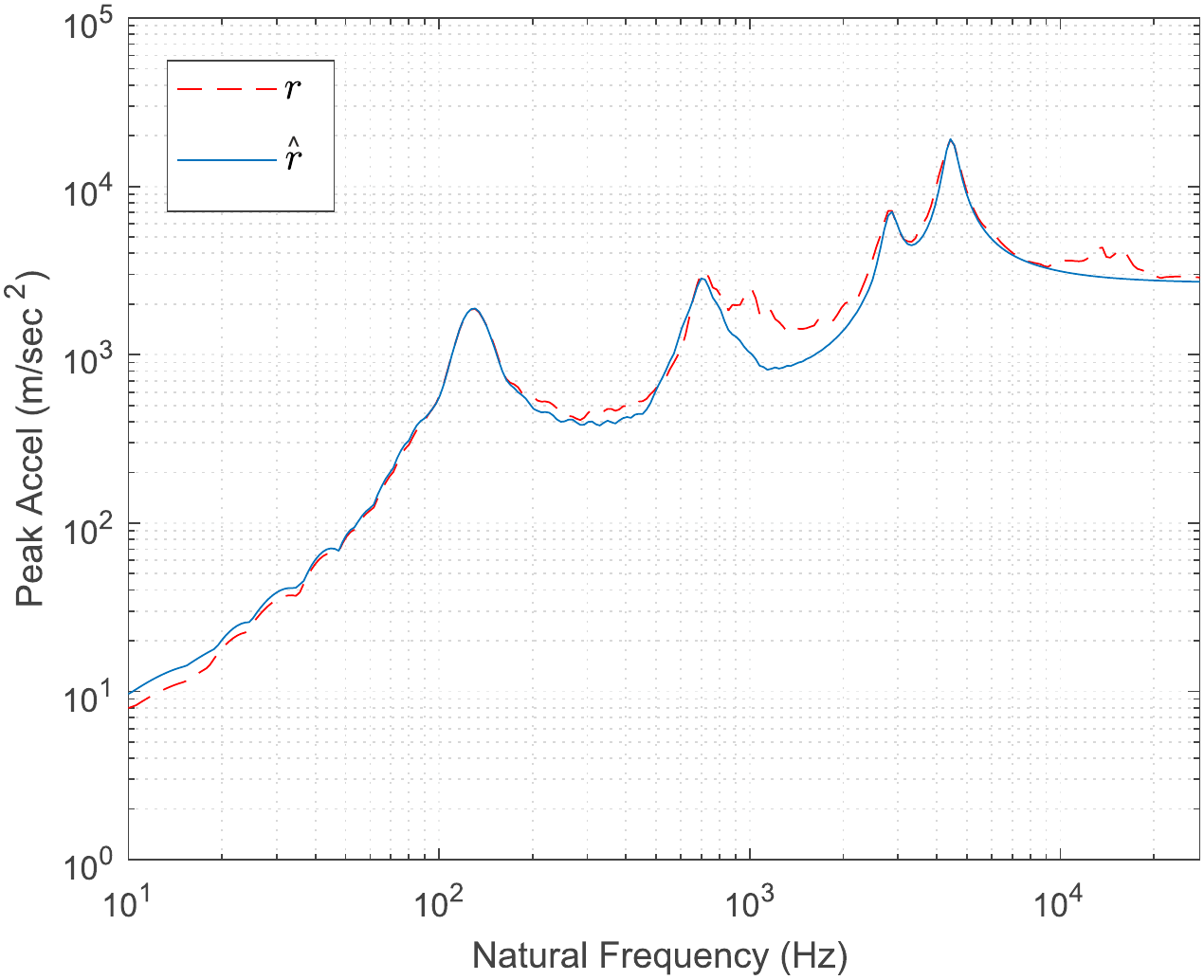}
		\caption{SRS}
		\label{srs_FFS}
	\end{subfigure}
	\caption{FFT and SRS of $r$ and $\hat{r}$ of FFS}
	\label{fft_srs_FFS}
\end{figure}

Validation of representativeness in time and frequency domain between simplified shock $\hat{r}$ and original shock signal $r$ are depicted in Figs.\ref{acc_FFS} and \ref{fft_srs_FFS}, respectively.
The signal $\hat{r}$ reconstructed from only 5 shock waveform components can keep both important temporal and frequency features of the original shock signal.

\begin{figure}
	\centering
	\includegraphics[width=\linewidth]{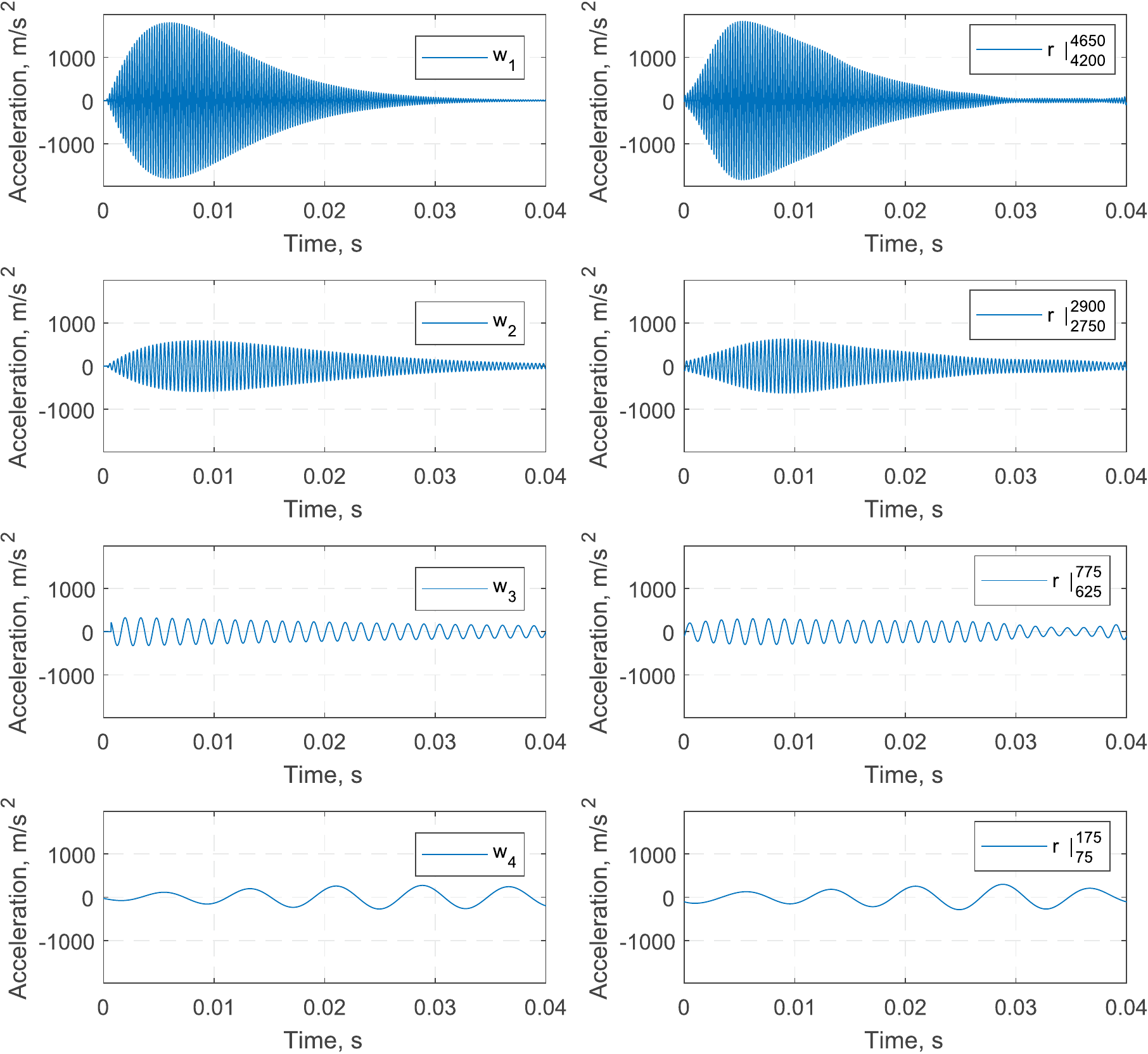}
	\caption{Comparison of $w_i$ and $r|_{\xi_i}^{\xi^i}$ for FFS}
	\label{bandpass_comparison}
\end{figure}

FFS is an appropriate example of Eq.(\ref{bandpass_relation}) for its frequency structure. 
The band of each bell $[\xi_i,\xi^i]$ can be easily identified in the FFT spectrum because the bells are clearly separated from each other ($\xi^i \leqslant \xi_{i+1}$), which are not satisfied in the case of NFS and MFS in this study.
Fig.\ref{bandpass_comparison} compares the decomposed shock waveform components and their corresponding filtered signals, which are extracted by the ideal band-pass filter in Matlab.
The pass-bands are shown in their respective legends.
It is evident that both Eqs.(\ref{bandpass_relation}) and (\ref{symmetry_relation}) are satisfied in this FFS case.
Band-pass filtered signals are almost the same as their corresponding shock waveform components.
Lower bound $\xi_i$ and upper bound $\xi^i$ of filtering band are symmetry about the centre frequency $w_i$.

\section{Conclusive Remarks}

According to the shock generation mechanism, this study proposes an explicit mathematical expression of the basic shock waveform, with which the time-frequency structure of a shock environment can be characterised, and the links between time and frequency domains can be identified.
Several index parameters with clear mechanical meaning are introduced, e.g., $\kappa$ and $\eta_{90\%}$, which can quantitatively describe the dominant shock distance and shock complexity.
By retaining only important waveform components, a shock environment can be simplified while still keeping its main mechanical properties in both time and frequency domain. 
By comparing specific parameters that are not easily available from other signal processing tools, the proposed method can be used to compare different shock environment and track the change of characteristics of the shock signal at different location when it propagates in a structure.
With the mathematical expression of shock signal, more convenient and reliable tool can be developed to characterise shock environment and its effects on equipment and devices.

\section*{Acknowledgements}
The authors would like to thank William Fitzpatrick, Dr. Stephen Burley and William Storey for the conduction and support of the experiment of mechanical impact.

\bibliographystyle{unsrt}  


\setcounter{figure}{0}
\appendix
\section{Mechanical Basis for the Proposed Shock Waveform}\label{mechanical_background}

Most of the shocks can be characterised by the free vibration of a structure subjected to a local impulsive loading.
The loading duration in most shock scenarios is short, comparing to the periods of the dominant vibration modes of the structure, which transmit the shock from the loading location to the concerned location.
This implies that the shock signal at the concerned location carries the contributions from the modes with high order natural frequencies of the structure.
If a structure is considered as a MDOF system, the shock response of such system is very complex, although the response of each modal is standard damped harmonic oscillation.
This is mainly because the modal density in high frequency domain (normally one order higher than the fundamental natural frequency) of target structure is very high, and the modal responses interact with each other.

\begin{figure}
	\centering
	\includegraphics[width=0.3\linewidth]{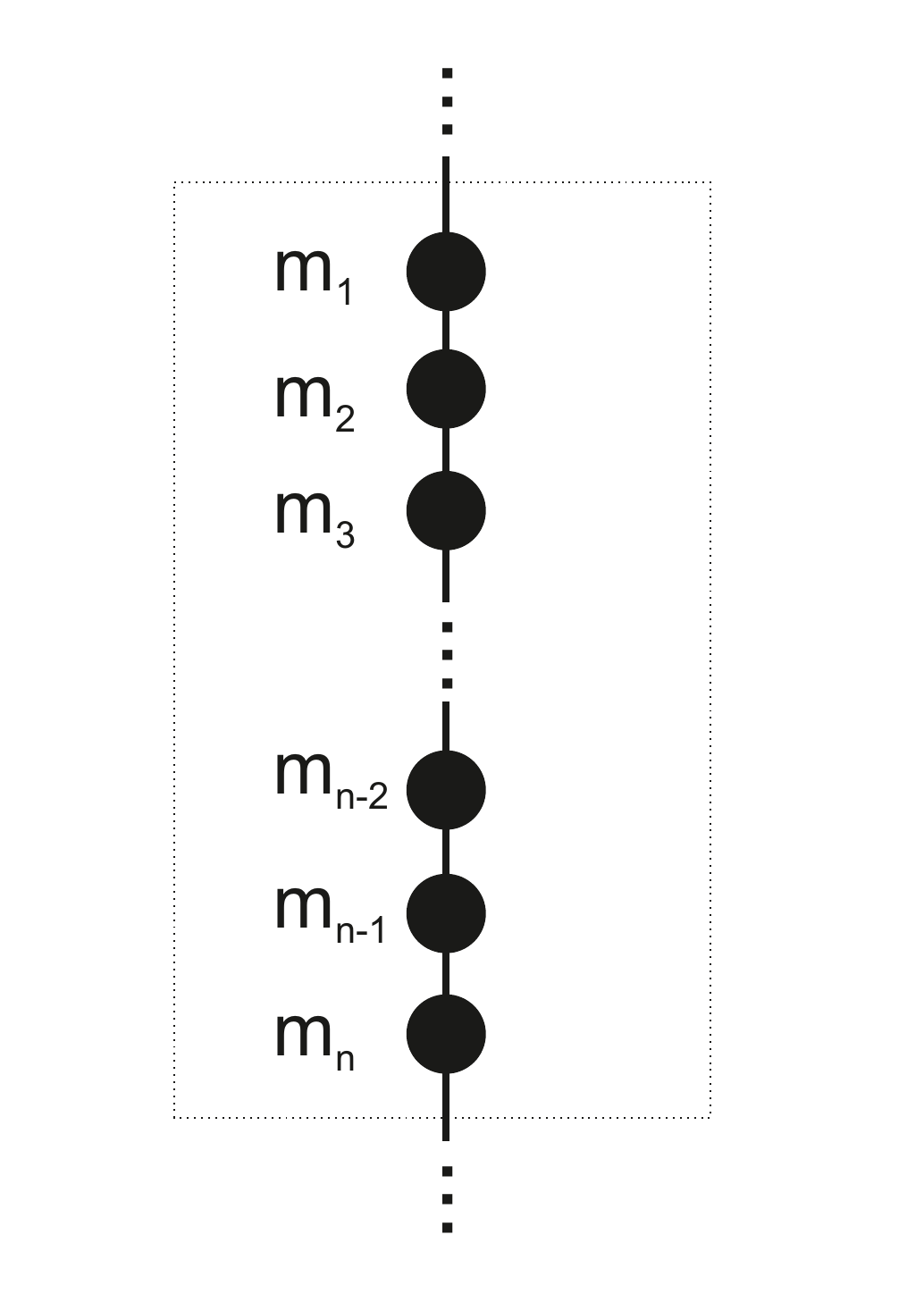}
	\caption{A multiple degree of freedom model (MDOF)}
	\label{ndof}
\end{figure}

It is assumed that a shock signal can be decomposed into a series of waveforms, and each of which can be represented by the response of a MDOF substructure subjected to the impulse loading.
The frequency band of each MDOF substructure is continuous and narrow.
Fig.\ref{ndof} shows a diagram of $n$-mass system representing a substructure.
As shown in Eq.(\ref{frequency_assumption}), the nature frequencies of the $j$th and $k$th modes are expressed as $\omega_j$ and $\omega_k$, respectively, which are assumed to be different but close to each other with a mean value $\bar{\omega}$.
The damping ratio of these $n$ modes are all assumed to be $\zeta$.
\begin{equation}\label{frequency_assumption}
\left\{
\begin{array}{lr}
\frac{1}{n}\sum_{j=1}^{n}\omega_j=\bar{\omega}\\
\frac{|\omega_j - \bar{\omega}|}{\bar{\omega}}\ll 1,\quad j \in [1,n]
\end{array}
\right.
\end{equation}

The dynamic motion equation of free vibration is
\begin{equation}
\m\ddot{\v}+\c\dot{\v}+\k\v=0
\label{motion_equation}
\end{equation}
where $\m$ is the mass matrix, $\c$ is the damping matrix, $\k$ is the stiffness matrix, and $\v$ is the displacement vector.

A solution in vector form for this motion equation is
\begin{equation}
\v=e^{\A t} \u_0
\end{equation}
where $\A$ is the solution matrix and $\u_0$ is the initial displacement vector.

The solution can be expressed with its eigenvector matrix $\Phi$.
\begin{equation}\label{complete_solution}
\v =\PHI e^{\LAMBDA_A t} \PHI^T  \u_0
\end{equation}
where $\LAMBDA_A$ is the eigenvalue matrix of matrix $A$.
According to the mechanical assumption in Eq.(\ref{frequency_assumption}), the diagonal elements of $\LAMBDA_A$ are $\lambda_j=-\omega_j\zeta+i\omega_j$, whose mean value is $\bar{\lambda}=-\bar{\omega}\zeta+i\bar{\omega}$.

Note that the eigenvalues are centralized around the mean $\bar{\lambda}$, the solution vector $\v$ can be written as
\begin{equation}
\v =e^{\A t} \u_0=e^{\bar{\lambda}\I t} e^{\B t} \u_0
\end{equation}
where
\begin{equation}
\B=\A-\bar{\lambda}\I
\end{equation}
According to the property of eigenvalue method, the eigenvalue matrix of $\B$ is
\begin{equation}
\LAMBDA_B=\LAMBDA_A-\bar{\lambda}\I
\end{equation}
Since the eigenvalue $\lambda_j$ in matrix $\LAMBDA_A$ is close to $\bar{\lambda}$, the eigenvalue $\lambda^B_j$ in matrix $\LAMBDA_B$ can be considered as small value.

The matrix $e^{\B t}$ can be expanded by Taylor series, i.e.
\begin{align}
e^{\B t} & =\sum_{n=0}^{\infty} \frac{1}{n!}(\B t)^n\label{taylor_series}\\
& = \I+\B t+\frac{1}{2}(\B t)^2+\frac{1}{6}(\B t)^3+\cdots\\
& = \I+\PHI\LAMBDA_B t\PHI^T+\frac{1}{2}\PHI(\LAMBDA_B t)^2\PHI^T+\frac{1}{6}\PHI(\LAMBDA_B t)^3\PHI^T+\cdots
\end{align}
The complete solution is
\begin{align}\label{final_solution}
\v & =\sum_{n=0}^{\infty}\v_n \nonumber\\
& =e^{\bar{\lambda}\I t} (\I+\PHI\LAMBDA_B t\PHI^T+\frac{1}{2}\PHI(\LAMBDA_B t)^2\PHI^T+\frac{1}{6}\PHI(\LAMBDA_B t)^3\PHI^T+\cdots) \u_0
\end{align}

The $t$ is relevant to the dominant frequency, i.e., $\bar{\omega}$ in this case.
In pyroshock, ballistic shock and navy shock scenarios, $t$ within shock duration is very small.
Combined with the small value of $\lambda^B_j$ in matrix $\LAMBDA_B$, the matrix series in Eq.(\ref{final_solution}) can converge fast.
By this way, the equations to describe such a shock response can be reduced form $n$ equations to much fewer equations.

For instance, if only the $n$th term in Eq.(\ref{taylor_series}) is outstanding, and only the response of the $m$th DOF is interested, then 
\begin{align}
\v_{n} & =\frac{1}{n!} e^{\bar{\lambda}\I t} (\B t)^n \u_{0}\\
& =\frac{t^n}{n!} e^{\bar{\lambda} \I t} \B^n \u_{0}\\
v_{mn} & =C t^n e^{\bar{\lambda} t}
\end{align}
where constant $C$ equals to the $m$th term in vector $\frac{1}{n!}\B^n \u_{0}$.

Recall that
\begin{equation}
\bar{\lambda}=-\bar{\omega}\zeta+i\bar{\omega}
\end{equation}
so the general waveform expression is in the form of
\begin{equation}
v_{mn} =C t^n e^{ -\bar{\omega}\zeta t+i\bar{\omega}t}
\end{equation}

\section{Discussion of Advanced Prony Mode}\label{discussion_AP_SW}

The advanced Prony mode is obtained by convoluting the standard Prony mode $P(t)$ with a basic Gaussian pulse $G(t)$ according to Refs.\cite{ECSS2015,bernaudin2005}. 
This section attempts to derive the explicit mathematical expression of advance Prony mode, which was not given in Refs.\cite{ECSS2015,bernaudin2005}. 

The following form of $P(t)$ and $G(t)$ are used
\begin{align}
P(t) & =A e^{-\zeta\omega t + i (\omega t + \varphi)}\\
G(t) & =e^{-\frac{(t-\tau )^2}{2\sigma ^2 \tau ^2}}
\end{align}
where $A$ is the amplitude, $\zeta$ is the damping ratio, $\omega$ is the angular frequency, $\varphi$ is the phase, $\tau$ is the peak time of Gaussian pulse, and $\sigma$ is a parameter relative to the width of Gaussian pulse.

The advanced Prony mode $P_{A} (t)$ can be obtain with Mathematica
\begin{align}
P_{A} (t) & =(P*G)(t)\\
& = \int_{0}^{t} G(x) P(t-x) dx\\
& =\int_{0}^{t} A \exp (-\zeta  \omega  (t-x)+i (\omega  (t-x)+\phi )-\frac{(x-\tau )^2}{2\sigma ^2 \tau ^2}) dx\\
&=A \sqrt{\frac{\pi }{2}} \sigma  \tau (\text{erf}(\frac{1+(\zeta -i) \sigma ^2 \tau  \omega }{\sqrt{2} \sigma })+\text{erf}(\frac{t-\tau  (1+(\zeta -i) \sigma ^2 \tau  \omega )}{\sqrt{2} \sigma  \tau })) \nonumber\\
&\quad \exp (\frac{1}{2} (\zeta -i) \omega  (-2 t+\tau  (2+(\zeta -i) \sigma ^2 \tau  \omega ))+i \phi )\label{expression_AP}
\end{align}
where $\text{erf}(x)$ is the Gauss error function.

The number of parameters in the explicit mathematical expression Eq.(\ref{expression_AP}) are the same as the proposed shock waveform in Eq.(\ref{shock_waveform}), if initial time $\mathring{t}$ is taken into account.
Although the expression of Eq.(\ref{expression_AP}) is very different from the proposed shock waveform in Eq.(\ref{shock_waveform}), the advanced Prony mode can also cover a wide range of existing waveform, e.g., the waveforms listed in Table \ref{special_cases}. 

\begin{figure}
	\centering
	\includegraphics[width=0.4\linewidth]{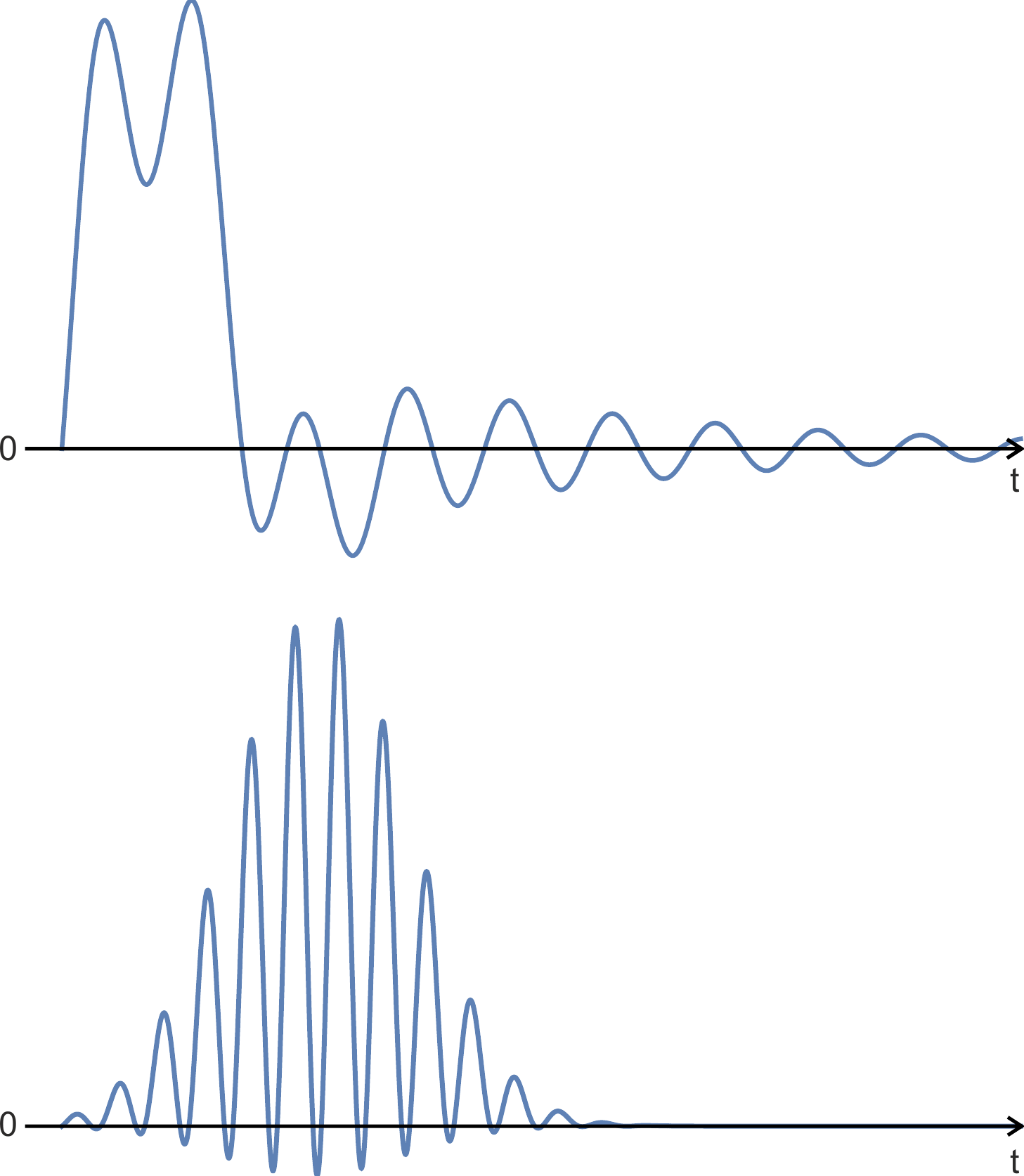}
	\caption{Examples of advanced Prony mode uncommonly found in mechanical shock environment}
	\label{abnormal_ap}
\end{figure}

However, the advanced Prony mode is not as convenient as the proposed shock waveform in Eq.(\ref{shock_waveform}) in terms of practical application.
It is difficult to get an intuitive understanding of mechanical meaning of each parameter.
The normalisation process of Eq.(\ref{expression_AP}) cannot be implemented, since the absolute maximum of this expression is difficult to be obtained.
The complexity mainly comes from the existence of time variable $t$ in the Gauss error function, which is a special (non-elementary) function.
Advanced Prony mode may even include some waveforms which are uncommon for mechanical shock.
Two such examples of waveforms produced by Eq.(\ref{expression_AP}) are shown in Fig.(\ref{abnormal_ap}) where the `zero velocity change' validation criteria\cite[p.417]{ECSS2015} is not satisfied.

\end{document}